\documentclass[12pt]{extarticle}
\usepackage[a4paper]{geometry}
\usepackage[cp1251]{inputenc}
\usepackage[russian,english]{babel}
\usepackage{StyleDeNiko}
\usepackage{amsbsy}
\usepackage{float}
\usepackage{graphicx}
\usepackage{threeparttable}
\usepackage{longtable}
\usepackage{booktabs,caption,fixltx2e,tabularx}
\usepackage{lscape}
\usepackage{multirow}
\usepackage{rotating}
\usepackage[sort&compress]{natbib}

\newcommand{\R}{\mathbb{R}}
\newcommand{\I}{\mathbb{I}}

\DeclareMathOperator{\inter}{int}

\newtheorem{Theorem}{Theorem}
\newtheorem{TheoremPrime}{Theorem}

\newtheorem{Coroll}{Corollary}
\newtheorem{CorollPrime}{Corollary}

\makeatletter
\renewcommand{\@biblabel}[1]{#1.}
\makeatother

\newcolumntype{L}[1]{>{\raggedright\arraybackslash}p{#1}}

\newcommand{\dpr}[2]{\left\langle #1,\,#2 \right\rangle}
\newcommand{\RP}{\mathit{R}^*(\omega)}
\newcommand{\ve}{\varepsilon}

\begin{document}

\begin{center}
  \LARGE{Positively-homogeneous Kon\"us-Divisia indices and their applications to demand analysis and forecasting}
\end{center}

\begin{center}
  Nikolay I. Klemashev

  \textit{Lomonosov Moscow State University}
\end{center}

\begin{center}
  and
\end{center}

\begin{center}
  Alexander A. Shananin

  \textit{Lomonosov Moscow State University}

  \textit{Moscow Institute of Physics and Technology}
\end{center}

{ \small \textit{Abstract.} This paper is devoted to revealed
preference theory and its applications to testing economic data for
consistency with utility maximization hypothesis, construction of
index numbers, and forecasting. The quantitative measures of
inconsistency of economic data with utility maximization behavior
are also discussed. The structure of the paper is based on
comparison between the two tests of revealed preference theory --
generalized axiom of revealed preference (GARP) and homothetic axiom
of revealed prefernce (HARP). We do this comparison both
theoretically and empirically. In particular we assess empirically
the power of these tests for consistency with maximization behavior
and the size of forecasting sets based on them. For the forecasting
problem we show that when using HARP there is an effective way of
building the forecasting set since this set is given by the solution
of the system of linear inequalities. The paper also touches upon
the question of testing a set of Engel curves rather than finite set
of observations for consistency with utility maximization behavior
and shows that this question has effective solution when we require
the rationalizing utility function to be positively homogeneous.}


{\small \textit{JEL Codes:} D11, D12, C14, C43

\section{Introduction}

Statistical services in various countries compute economic indices
using the approach based on the value of consumption bundle in
different time periods. The problem of this approach is that the
consumption bundle changes with the change of the structure of
prices. The statistical services compute Laspeyres and Paasche price
indices\footnote{See \citep{Diewert2004ILOCh15}, for example} fixing
some consumption bundle equal to consumption levels in base or
current period. More precisely, for Laspeyres index they take the
set of goods consumers bought in base period as consumption bundle.
For Paasche index consumption in current period is taken as the
consumption bundle. Usually the value of Paasche index is not
greater than that of Laspeyres index. This systematic difference is
called Gerschenkron's effect. It is caused by substitution effects
-- when a price on certain commodity increases people tend to
substitute it with different goods. This implies that consumption
structure is changing and no fixed consumption bundle can represent
what actually consumers buy. See \citep{CPIManual2004,
Ershov2011Eng} for a detailed review of existing approaches to the
problem of index numbers construction.

We address the issue of substitution effects by suggesting a new
approach to index numbers. The approach is based on Pareto's theory
of demand. According to this theory a representative economic agent
is supposed to choose his or her consumption bundle among
inter-substitutable goods by maximizing his or her utility function
given his or her budget constraint. Let us call the problem of
construction demand functions from utility function as a direct
problem. Then the inverse one is to recover utility function from
demand functions. Solvability of inverse problem is related to
integrability of demand functions. This problem is central for
revealed preference theory.

An approach for construction of consumption index based on the
solution of inverse problem was suggested in
\citep{Konus1924Eng}\footnote{See \citep{Konus1939} for English
translation}. This index is called Kon\"us consumption index. An
effective way of computing Kon\"us consumption index is based on
Afriat's contributions to revealed preference theory (see
\citep{Afriat1967}).

Divisia index numbers formulas, first introduced in
\citep{Divisia1925}, represent a general form of consumption and
price indices. See \citep{Balk2005} for further details, systematic
description and a historical overview of Divisia index numbers. In
general, Kon\"us consumption index is not a Divisia index. However,
when we add the requirement for the utility function to be
positive-homogeneous of degree~1, Kon\"us consumption index becomes
Divisia index \citep{Hulten1973} and is called Kon\"us-Divisia
consumption index. Moreover, one may construct a dual index which is
interpreted as price index. This index is called Kon\"us-Divisia
price index. An interesting property of Kon\"us-Divisia price and
consumption indices is that they are related through Yang transform
\citep{Shananin1986Eng}.

Revealed preference theory allows to test whether a finite set of
observations on prices and consumption for a certain group of goods
(this is called trade statistics) is consistent with Pareto demand
theory. The criterion for the trade statistics to be consistent with
Pareto theory with nonnegative, nonsatiated, continuous, concave,
and monotonic utility function is generalized axiom of revealed
preference (GARP). The criterion for the trade statistics to be
consistent with Pareto theory with additional requirement for the
utility function to be positively-homogeneous of degree~1 is
homothetic axiom of revealed preference (HARP). If the trade
statistics is consistent with Pareto demand theory then we say that
it (or the corresponding group of goods) is rationalizable by a
utility function. For treatment of the model of collective
consumption when one explicitly model total household consumption as
the result of decisions made by its members see
\citep{CherchyeDeRockVermeulen2011, CherchyeDeRockVermeulen2012,
CherchyeDeRockPlatino2013}. For developments in revealed preference
techniques for applications other than rationalizability of consumer
behavior on consumption commodities markets see \citep{Varian2012,
EcheniqueEtAl2013, CarvajalEtAl2013}.

Rationalizability of the group means that it is full in terms of
substitutes and complements. For any good in such group all the
goods which are substitutes or complements to it should also belong
to the group. The concept of rationalizabilty is related to the
concept of separability. A subgroup of a rationalizable group is
said to be separable if it is rationalizable. The separability
allows one to build a structure of consumer demand.

The rationalizability of a group of goods should have some intuitive
explanation. In other words one should be able to explain why a
rationalizable group is full in terms of substitutes and
complements. If almost any group of goods is rationalizable then the
concept of rationalizability does not provide a tool for studying
the segmentation of market, because it allows almost any structure
of consumer demand. We show in this paper that this is the case for
rationalizability based on GARP. More precisely, we show empirically
that a randomly chosen group of goods satisfies GARP with high
probability. Several researchers (see \citep{BeattyCrawford2011},
for example) notice that one of the problem of low power of GARP is
that the prices and consumption data are growing over time due to
economic growth and this implies low number of budget sets
intersections. This is the reason why in practice we often are
unable to reject GARP.

The situation is different with HARP because this axiom is invariant
with respect to changing of scales in economic data and the low
number of intersections between budget sets does not lower the
restrictiveness of HARP. We provide empirical evidence in favor of
this conclusion by showing that the probability of having a random
group (of not very large size) satisfying HARP is very low. This
implies that if a group satisfies HARP then most probably it is not
a mere coincidence and the group indeed may be interpreted as an
independent market segment.

%Trade statistics is a set of observations on consumption and prices
%for a certain group of goods during a fixed time interval. That is
%why one may misclassify non-rationalizable groups as rationalizable
%ones because of not having enough observations. This implies that
%the constructed structure of consumer demand may be subject to
%serious changes as we expand the time interval and add more
%observations on consumption and prices. In order to have a stable
%structure of consumer demand we need to select groups which remains
%rationalizable as we add new observations. One of directions of
%changing prices and consumption is changing of scales.
%
These observations led us to conjecture that one should use HARP as
the instrument of revealing the true market segmentation. This means
that one should consider rationalizability only by
positively-homogeneous utility functions. In further sections we
provide more empirical evidence in favor of importance of the
requirement for the utility function to be positively-homogeneous by
studying the power of tests based both on GARP and HARP.

If a trade statistics does not satisfy GARP and/or HARP there may be
two different reasons for this. The first reason is that the
selected group of goods is not complete and the consumer behavior
cannot be explained by means of the concept of rational
representative consumer. The second possible cause is the
measurement error in data. In order to differ these two reasons one
needs to have a quantitative measure of the degree to which the
trade statistics satisfies revealed preference axioms. The first
such measure was suggested by \citep{Afriat1973}. After that there
were many papers devoted to alternative ways of measuring the degree
of violation of trade statistics with GARP (see
\citep{WhitneySwofford1987, Varian1990, Famulari1995, Gross1995,
EcheniqueLeeShum2011, SmeuldersEtAll2013, EkelandGalichon2013}).

The power of GARP was assessed in \citep{Bronars1987}. The author
generated random consumption data and then checked how often the
data set constructed from original prices and randomized consumption
satisfied GARP. The drawback of the author's approach is that he
ignores correlation in time series of consumption. He generated
random consumption data independently for each period. We eliminate
this drawback by using random process in generating randomized data
sets. We have not seen empirical papers on the power of HARP. We
also have not seen any study of the probability of having a random
group satisfying HARP and GARP.

The two tests in revealed preference theory allow one to predict
consumption at some fixed price. The approach which is due to
\citep{Varian1982a} is to find all possible consumption bundles
which when combined with the given data set constitute
rationalizable trade statistics. The similar approach is used to
predict prices given some consumption bundle. In applied research
one prefers to have such set of possible price vectors to be as
small as possible. We provide empirical evidence showing that the
approach based on GARP may result in very large sets which are close
to the set of all vectors with positive elements. However, the
approach based on HARP provide much smaller sets for the same groups
of commodities.

We measure the size of the set of predicted prices as the
probability of having a trade statistics with the last period price
being a random vector to satisfy revealed preference axiom (GARP or
HARP). This approach is close to the measure of restrictiveness of
GARP suggested by \citep{BeattyCrawford2011} who also suggest the
measure of success in testing GARP which takes into account the
strictness of the constraints implied by GARP. They also report low
level of restrictiveness of GARP constraints from their empirical
study. These findings are consistent with our findings about the
size of the set of predicted prices implied by GARP.

Having a set of Engel curves we say that this set is rationalizable
in some class of utility functions if any trade statistics
consisting of the prices defining the Engel curves and set of
demands on these curves (one demand from each curve) is
rationalizable in the same class of utility functions. The problem
of testing a set of arbitrary Engel curves for consistency with GARP
is still an open one. However, in case when Engel curves are rays
this problem has solution. Namely, in order to test a set such Engel
curves for consistency with GARP one needs to test a single trade
statistics built from these curves for consistency with HARP. For
arbitrary set of Engel curves, even without intersections and
monotone, it is not enough to check a single trade statistics for
consistency with GARP to prove their consistency with GARP.

One may use the concept of rationalizability of Engel curves in
order to build a smaller forecasting set for prices or consumption.
In case of HARP, however, this does not lead to narrowing of the
forecasting set. This may make forecasting set based on consistency
with GARP smaller. However, since there is no way of checking for
consistency of Engel curves with GARP\footnote{Except for checking
all possible trade statistics built from the set Engel curves for
consistency with GARP.}, there is no way to construct this
forecasting set even numerically. An attempt to provide an algorithm
for building this set was made by
\citep{BlundellBrowningCrawford2008, BlundellKristensenMatzkin2014}.
However, the method suggested there is incorrect. We formulate the
main result of \citep{BlundellBrowningCrawford2008} in
subsection~\ref{SSForecasting} and prove its incorrectness by
providing a counter-example in appendix~2.

The rest of the paper is organized in the following way. In section
\ref{STheory} we describe the two tests for consistency of data with
utility maximization hypothesis, their generalizations on the case
of data with measurement errors, how to forecast prices or
consumption, and the notions of separability and hierarchy of index
numbers\footnote{This also is known as utility tree.}. Section
\ref{SApplics} contains examples of applications of Kon\"us-Divisia
index numbers to real data. In section \ref{SPower} we estimate the
power of the tests, and in section \ref{SSize} we assess the sizes
of forecasting sets based on the considered tests from revealed
preference theory. Section \ref{SConcl} concludes. Appendix~1
contains proofs for yet unpublished results. Appendix~2 provides
counter-example showing incorrectness of theorem~1 in
\citep{BlundellBrowningCrawford2008}.

\section{Theory}\label{STheory}

In this section we describe the instruments for testing data for
consistency with maximization behavior and forecasting. We also
describe a quantitative measures of the degree to what a data set
fails to satisfy the hypothesis of maximization behavior. The
methods of forecasting and building a tree of economic indices are
also discussed.

\subsection{Testing data for consistency with utility maximization
hypothesis} \label{HARPGARP}

Consider a group of $m$ goods. Denote a consumption bundle of these
goods by $X = (X_1,\ldots,X_m)$ and their prices by $P =
(P_1,\ldots,P_m)$. A finite set $\{(P^t, X^t)\}_{t=1}^T$ of
observations on prices and consumption in periods $1,\ldots, T$ is
called trade statistics. We say that it is rationalizable in
functional class $\Phi$ if there exists a utility function $u \in
\Phi$ such that
\begin{gather*}
  X^t\in\Argmax\left\{u(X)\left| \dpr{P^t}{X}\leqslant \dpr{P^t}{X^t}, X\in\R_+^m\right.\right\} \quad \forall t \in \{1,\ldots,T\},
\end{gather*}
where $\dpr{\cdot}{\cdot}$ denotes scalar product and $\R_+^m =
\{X\in\R^m \mid X \geqslant 0, X \ne 0\}$.

In this paper we consider two classes of utility functions. The
first one is that of nonnegative, nonsatiated, continuous, concave,
and monotonic on $\R^m_+$ and positive on $\inter\R_+^m$ functions.
We denote it by~$\Phi_G$. The second one is that of functions
from~$\Phi_G$ which are positively-homogeneous of degree 1. The
second class is denoted as~$\Phi_H$.

The trade statistics $\{(P^t,X^t)\}_{t=1}^T$ may be tested for
rationalizability in these two classes by means of the two
Afriat-Varian theorems.

\begin{Theorem}{\citep{Afriat1963, Afriat1967, Varian1983}}\label{AVGARPTh}
  The following statements are equivalent
  \begin{itemize}
    \item[1)] trade statistics $\{(P^t, X^t)\}_{t=1}^T$ is rationalizable in $\Phi_G$;
    \item[2)] there exist numbers $U^t$, $\lambda^t > 0$,
    $(t=1,\ldots,T)$, such that
    \begin{gather}
      U^t \leqslant U^s + \lambda^s\left(\dpr{P^s}{X^t} - \dpr{P^s}{X^s}\right), \quad \forall
      t,s = 1,\ldots,T; \label{GARPSys}
    \end{gather}
    \item[3)] trade statistics $\{(P^t, X^t)\}_{t=1}^T$ satisfies GARP, which means that for all $t$
    and $s$ if there exists $t_1,\ldots,t_k$ such that
    $\dpr{P^t}{X^t}\geqslant \dpr{P^t}{X^{t_1}}$, $\dpr{P^{t_1}}{X^{t_1}}\geqslant
    \dpr{P^{t_1}}{X^{t_2}}$, $\ldots$, $\dpr{P^{t_k}}{X^{t_k}}\geqslant \dpr{P^{t_k}}{X^s}$,
    then $\dpr{P^s}{X^s} \leqslant \dpr{P^s}{X^t}$;
    \item[4)] the function $F_G(X) = \min_{s\in\{1,\ldots,T\}}\left\{U^s + \lambda^s\left(\dpr{P^s}{X} -
    \dpr{P^s}{X^s}\right)\right\}$, where $\{U^s, \lambda^s\}_{s=1}^T$ satisfy (\ref{GARPSys}) and
    $\lambda^s > 0$ $\forall s \in \{1,\ldots,T\}$, rationalizes trade statistics.
  \end{itemize}
\end{Theorem}

In order to check whether the trade statistics satisfies GARP one
needs to build a transitive closure of the relation $\mathit{R}$
defined on $\{X^t\}_{t=1}^T\times\{X^t\}_{t=1}^T$ as
\begin{gather}
  X^t\mathit{R}X^s \Leftrightarrow \dpr{P^t}{X^t} \geqslant \dpr{P^t}{X^s}.
  \label{RevealedPref}
\end{gather}
If the resulting transitive closure $\mathit{R}^*$ is such that for
all $t, s \in \{1,\ldots,T\}$ $X^tR^*X^s$ implies $\dpr{P^s}{X^s}
\leqslant \dpr{P^s}{X^t}$, then the trade statistics satisfies GARP.
See \citep{Varian1982a} for an algorithm which allows one to
construct a solution to the system of linear inequalities
(\ref{GARPSys}) if the trade statistics satisfies GARP.

\begin{Theorem}{\citep{Afriat1963, Afriat1967, Varian1983}}\label{HARPFinTh}
The following statements are equivalent
  \begin{itemize}
    \item[1)] trade statistics $\{(P^t, X^t)\}_{t=1}^T$ is rationalizable in $\Phi_H$;
    \item[2)] there exist numbers $\lambda^t > 0$ $(t =
    1,\ldots,T)$, such that
    \begin{gather}
      \lambda^t\dpr{P^t}{X^s} \geqslant \lambda^s\dpr{P^s}{X^s}, \quad \forall t,s =
      1,\ldots,T; \label{HARPSys}
    \end{gather}
    \item[3)] trade statistics $\{(P^t, X^t)\}_{t=1}^T$ satisfies Homothetic Axiom of Revealed
    Preference (HARP), which means that for all subsets of indices
    $\{t_1,\ldots,t_k\}$ from $\{1,\ldots,T\}$ the following
    inequality is satisfied:
    \begin{multline}
      \dpr{P^{t_1}}{X^{t_2}} \dpr{P^{t_2}}{X^{t_3}} \ldots
      \dpr{P^{t_k}}{X^{t_1}} \\ \geqslant \dpr{P^{t_1}}{X^{t_1}} \dpr{P^{t_2}}{X^{t_2}} \ldots
      \dpr{P^{t_k}}{X^{t_k}}. \label{HARPIneq}
    \end{multline}
    \item[4)] the function $F(X) =
    \min_{s\in\{1,\ldots,T\}}\left\{\lambda^s\dpr{P^s}{X}\right\}$, where
    $\{\lambda^s\}_{s=1}^T$ satisfy (\ref{HARPSys}) and $\lambda^s >
    0$ $\forall s \in \{1,\ldots,T\}$,
    rationalizes trade statistics.
  \end{itemize}
\end{Theorem}

The inequality (\ref{HARPIneq}) for $k = 2$ may be written as
\begin{gather*}
  \frac{\dpr{P^{t_1}}{X^{t_2}}}{\dpr{P^{t_2}}{X^{t_2}}} \geqslant
  \frac{\dpr{P^{t_1}}{X^{t_1}}}{\dpr{P^{t_2}}{X^{t_1}}}.
\end{gather*}
This inequality reflects Gerschenkron's effect -- Laspeyres index is
no less than Paasche index ($t_2$ is base period, $t_1$ is current
period).

A solution to the system of linear inequalities (\ref{HARPSys}) may
be found by means of Floyd-Warshall algorithm \citep{Warshall1962,
Floyd1962}. Introduce Paasche price indices matrix $C$ with elements
given by
\begin{gather*}
  C_{ts} = \frac{\dpr{P^s}{X^s}}{\dpr{P^t}{X^s}}.
\end{gather*}
Then the system (\ref{HARPSys}) may be rewritten as
\begin{gather}
  \lambda^sC_{ts} \leqslant \lambda^t. \label{HARPSysMod}
\end{gather}
Let $C^*_{ts}$ be equal to the maximum value of
$C_{tt_1}C_{t_1t_2}\ldots C_{t_ks}$ over all possible ordered
subsets $\{t_1,\ldots,t_k\}$ of $\{1,\ldots,T\}$ for any positive
integer $k$ assuming that the empty subset corresponds to $C_{ts}$,
that is
\begin{gather*}
  C_{ts}^* = \max\{C_{tt_1}C_{t_1t_2}\ldots C_{t_ks} \mid \{t_1,\ldots,t_k\}
  \subset \{1,\ldots,T\}, k \geqslant 0\}.
\end{gather*}
Trade statistics satisfies HARP if and only if $C_{tt}^* \leqslant
1$ $\forall t \in \{1,\ldots,T\}$. One may notice that if the system
(\ref{HARPSysMod}) has positive solution, then it is equivalent to
\begin{gather}
  \lambda^sC_{ts}^* \leqslant \lambda^t. \label{l16}
\end{gather}

Consider an idempotent semi-ring with operations $a\oplus b =
\max(a,b)$ and $a\otimes b = ab$. Then the matrix $C^*$ with the
elements $C_{ts}^*$ is given by
\begin{gather}
  C^* = C \oplus C^{\bullet2} \oplus \ldots \oplus C^{\bullet k}
  \oplus \ldots, \label{CStar}
\end{gather}
where $C^{\bullet n}$ means taking $n$-th power of the matrix $C$ in
idempotent sense with all summing operations replaced by $\oplus$.
Notice that all the elements of the matrix $C$ are positive. If at
some step $n \leqslant T$ of taking idempotent powers we have
$C_{tt}^*
> 1$ for some $t$, then all the elements of $C^*$ are equal to
$\infty$ and the system (\ref{l16}) has no solution. Otherwise, it
is enough to take the first $T$ terms in (\ref{CStar}) in order to
compute $C^*$. Therefore, the complexity of algorithm for
computation of $C^*$ is of order $T^3$. If the system (\ref{l16})
has a positive solution then
\begin{gather*}
  \lambda^t = \max\{C^*_{t\tau} \mid \tau \in \{1,\ldots,T\}\} \quad
  t = 1,\ldots,T
\end{gather*}
solve (\ref{HARPSysMod}), which means that these $\lambda^t$ solve
(\ref{HARPSys}).

Given a positive solution of the system (\ref{HARPSys}) one may
define Kon\"us consumption index as $F(X) =
\min\left\{\lambda^s\dpr{P^s}{X} \mid s\in\{1,\ldots,T\}\right\}$.
This index is also Divisia index. That is why we call it as
Kon\"us-Divisia consumption index. The corresponding Kon\"us-Divisia
price index is given by function $Q(P)$ which is Yang transform of
$F(X)$, that is
\begin{gather*}
  Q(P) = \inf\left\{\left.\frac{\dpr{P}{X}}{F(X)}\right|X\geqslant0, F(X) >
  0\right\}.
\end{gather*}
This transform comes from the optimization problem
\begin{gather*}
  \max\left\{F(Y)\left| \dpr{P(X)}{Y} \leqslant \dpr{P(X)}{X},\, Y \in \R_+^m\right.\right\},
\end{gather*}
where Lagrange multiplier to the constraint $\dpr{P(X)}{Y} \leqslant
\dpr{P(X)}{X}$ is given by $\frac1{Q(P(X))}$. Yang transform is
involutive in $\Phi_H$ (see \citep{Kleiner1980Eng}). In
\citep{VratenkovShananin1991Eng} it was suggested to compute
consumption and price Kon\"us-Divisia index numbers as $F(X^t) =
\lambda^t\dpr{P^t}{X^t}$ and $Q(P^t) = \frac1{\lambda^t}$.

This approach to construction of index numbers is called
nonparametric one (see \citep{Diewert1973, Shananin1993Eng,
Houtman1995}).

%In order to prove that the value of $Q(P)$ at $P^t$ is given by
%$\frac1{\lambda^t}$ one need to start with inequality
%\begin{gather*}
%  H(X) = \frac{P^tX}{\min\{\lambda^{\tau}P^{\tau}X\mid \tau \in
%  \{1,\ldots,T\}\}} \geqslant \frac{P^tX}{\lambda^tP^tX} =
%  \frac1{\lambda^t},
%\end{gather*}
%and notice that $H(X^t) = \frac1{\lambda^t}$.
%
%The Yang transform has the following properties (see
%\citep{Ashmanov1984Eng}):
%\begin{itemize}
%  \item[1)] The function $Q(P)$ is defined for any $P\in\R^m_+$.
%  \item[2)] $Q(P) \geqslant 0$.
%  \item[3)] $Q(P)$ is positively-homogeneous.
%  \item[4)] $Q(P)$ is concave.
%  \item[5)] $Q(P)F(X) \leqslant PX, \quad \forall X\in\R^m_+,
%  P\in\R_+^m$.
%  \item[6)] $Q(P)$ is not equal to zero for some $P\in\R_+^m$.
%  \item[7)] $Q(P)$ is continuous on $\R_+^m$.
%  \item[8)] Let $F^{\times\times}(X)$ be Yang transform of $Q(P)$, that
%  is
%  \begin{gather*}
%    F^{\times\times}(X) = \inf\left\{\left.\frac{PX}{Q(P)} \right| P \geqslant
%    0, Q(P) > 0\right\}.
%  \end{gather*}
%  If $F \in \Phi_H$ then $F^{\times\times}(X) = F(X)$.
%\end{itemize}
%The properties 1-4, 6, and 7 imply that Yang transform
%transforms functions from $\Phi_H$ to functions from $\Phi_H$. This
%property combined with the property 8 means that Yang transform
%is involutive in $\Phi_H$.

The relation between GARP and HARP is given by
\begin{Theorem}{\citep{Shananin2009Eng}}\label{HARPandGARP}
  The following statements are equivalent:
  \begin{itemize}
    \item[1)] trade statistics $\{(P^t, X^t)\}_{t=1}^T$ satisfies HARP;
    \item[2)] trade statistics $\{(P^t, \mu^tX^t)\}_{t=1}^T$ satisfies
    GARP for any positive values of $\mu^t$ $(t \in \{1,\ldots,T\})$.
  \end{itemize}
\end{Theorem}

In order to discuss one of the implications of this theorem we need
to define the concept of rationalizability of the set of Engel
curves. Consider a set of Engel curves $\{q^t(\cdot)\}_{t=1}^T$
corresponding to price vectors $\{P^t\}_{t=1}^T$. This set of Engel
curves is said to be rationalizable in $\Phi$ if there exists a
function $u\in\Phi$ rationalizing trade statistics $\{(P^t,
q^t(x_t))\}_{t=1}^T$ for any positive total expenditures $x_1$,
$\ldots$, $x_T$. Theorem~\ref{HARPandGARP} implies that if a trade
statistics $\{(P^t, X^t)\}_{t=1}^T$ satisfies HARP then the set of
Engel curves $\{q^t(\cdot)\}_{t=1}^T$ corresponding to prices
$\{P^t\}_{t=1}^T$ is rationalizable in $\Phi_G$. The fact that it is
also rationalizable in $\Phi_H$ comes from the invariance of HARP
with respect to change of scales in consumption data. Indeed, if the
trade statistics satisfies HARP then it is rationalizable by
positively-homogeneous function from $\Phi_G$. This implies that all
Engel curves are rays. Therefore, choosing trade statistics
corresponding to a different values of expenditures from the set of
Engel curves is equivalent to the change of scales in consumption
data.

The invariance of HARP with respect to change of scales in prices
and consumption data is a very important property. It makes
rationalizability in $\Phi_H$ a more reliable property since trade
statistics data are subject to changes of scales in time caused by
economic growth and inflation.

\subsection{Testing inverse demand functions for consistency with utility maximization
hypothesis and the hierarchy of Kon\"us-Divisia
indices}\label{HARPGARPInf}

In this subsection we provide modified versions of theorems
\ref{AVGARPTh} and \ref{HARPFinTh} for the case when we know the
inverse demand functions $P(X)$ rather than their values at finite
set of consumption bundles. We assume that the inverse demand
functions are continuous on $\R^m_+$.

The relation $\mathit{R}$ given by (\ref{RevealedPref}) is called a
direct revealed preference relation. It was introduced in
\citep{Samuelson1938} in case when we have inverse demand function
$P(X)$ as
\begin{gather*}
  X\mathit{R}Y \Leftrightarrow \dpr{P(X)}{X} \geqslant \dpr{P(X)}{Y}.
\end{gather*}

We say that the inverse demand functions are rationalizable in
functional class $\Phi$ if there exists a utility function $u \in
\Phi$ such that
\begin{gather*}
  X \in \Argmax\left\{u(Y)\left| \dpr{P(X)}{Y} \leqslant \dpr{P(X)}{X}, Y \in
  \R^m_+\right.\right\}, \forall X \in R_+^m.
\end{gather*}

\begin{Theorem}{\citep{Stigum1973}}
  Let $P(X)$ be nonnegative, continuous on $\R^m_+$ functions such
  that $\dpr{P(X)}{X} > 0$ for all $X \in \R_+^m$. Then the
  following statements are equivalent
  \begin{itemize}
    \item[1)] the inverse demand functions are rationalizable in
    $\Phi_G$;
    \item[2)] there exist continuous on $\mathop{int}\R^m_+$ functions
    $U(X)\in\Phi_G$ and $\lambda(X)$ with $\lambda(X) > 0$ $\forall X \in
    \R^m_+$, such that
    \begin{gather}
      U(X) \leqslant U(Y) + \lambda(Y)\left(\dpr{P(Y)}{X} - \dpr{P(Y)}{Y}\right) \quad \forall
      X \in \R_+^m, Y \in \R_+^m; \label{GARPSysInf}
    \end{gather}
    \item[3)] the inverse demand functions satisfy GARP, that is if for some
    $X$, $X^1$, $\ldots$, $X^K$, $Y$ from $\R_+^m$
    $\dpr{P(X)}{X}$ $\geqslant$ $\dpr{P(X)}{X^1}$, $\dpr{P(X^1)}{X^1} \geqslant
    \dpr{P(X^1)}{X^2}$, $\ldots$, $\dpr{P(X^K)}{X^K} \geqslant \dpr{P(X^K)}{Y}$ then
    $\dpr{P(Y)}{Y}
    \leqslant
    \dpr{P(Y)}{X}$.
  \end{itemize}
\end{Theorem}

Consider a function $F$ from $\Phi_H$. Let $\Hat P \in \partial
F(\Hat X)$, where $\partial F(\Hat X)$ is the superdifferential of
the function $F$ at point $\Hat X$. Then as shown in
\citep{Shananin2009Eng} the generalized Euler identity
\begin{gather*}
  F(\Hat X) = \dpr{\Hat P}{\Hat X}
\end{gather*}
holds true.

Let us look at the inequalities (\ref{GARPSysInf}). They imply that
$Y$ solves the optimization problem
\begin{align*}
  & U(X) \rightarrow \max\limits_{X}, \\
  & \text{s.t. } \dpr{P(Y)}{X} \leqslant \dpr{P(Y)}{Y}, \\
  & \phantom{\text{s.t. }} X \in \R_+^m.
\end{align*}
The first order condition for this problem is given by
\begin{gather*}
  \lambda(Y)P(Y) \in \partial U(Y).
\end{gather*}
If we require $U\in\Phi_H$ then the generalized Euler identity
implies
\begin{gather*}
  U(Y) = \lambda(Y)\dpr{P(Y)}{Y}  \quad \forall Y \in \R_+^m.
\end{gather*}
If we substitute this in (\ref{GARPSysInf}) we obtain the following
inequalities
\begin{gather*}
  \lambda(Y)\dpr{P(Y)}{X} \geqslant \lambda(X)\dpr{P(X)}{X}.
\end{gather*}
These observations provide intuition on the criteria for
rationalizability in $\Phi_H$.

\begin{Theorem}{\citep{Levin1997, PospelovaShananin1998Eng}}
\label{HARPThInf}
  Let $P(X)$ be nonnegative, continuous on $\R^m_+$ functions such
  that $\dpr{P(X)}{X} > 0$ for all $X \in \R_+^m$. Then the
  following statements are equivalent
  \begin{itemize}
    \item[1)] the inverse demand functions $P(X)$ are rationalizable in
    $\Phi_H$;
    \item[2)] the system of linear inequalities
    \begin{gather}
      \lambda(Y)\dpr{P(Y)}{X} \geqslant \lambda(X)\dpr{P(X)}{X} \quad (X\in\R_+^m, Y\in\R_+^m) \label{HARPSysInf}
    \end{gather}
    has a solution $\lambda(X)$ which is positive and continuous on
    $\mathop{int}\R^m_+$;
    \item[3)] the inverse demand functions $P(X)$ satisfy HARP, that
    is for any set of vectors $\{X^1, \ldots, X^T\}$ from $\R^m_+$
    the inequality
    \begin{multline*}
      \dpr{P(X^1)}{X^2} \dpr{P(X^2)}{X^3} \ldots \dpr{P(X^T)}{X^1}
      \\
      \geqslant \dpr{P(X^1)}{X^1} \dpr{P(X^2)}{X^2} \ldots
      \dpr{P(X^T)}{X^T}
    \end{multline*}
    holds true;
    \item[4)] there exist price index function $Q(P)$ and
    consumption
    index function $F(X)$ both from $\Phi_H$ such that
    \begin{align}
      & Q(P)F(X) \leqslant \dpr{P}{X} & & \forall P\in\R^m_+,
      X\in\R_+^m,  \label{HARPIndIneq1} \\
      & Q(P(X))F(X) = \dpr{P(X)}{X} & & \forall X \in \R_+^m. \label{HARPIndIneq2}
    \end{align}
  \end{itemize}
\end{Theorem}

\begin{Coroll} \label{HARPInfCol1}
    If $Q(P)$ and $F(X)$ are functions from $\Phi_H$ satisfying
    (\ref{HARPIndIneq1}) and (\ref{HARPIndIneq2}) then $\lambda(X) =
    \frac1{Q(P(X))}$ satisfies (\ref{HARPSysInf}).
\end{Coroll}

\begin{Coroll} \label{HARPInfCol2}
  The functions $Q(P(X))$ and $F(X)$ in 4) may be expressed through a
  solution $\lambda(X)$ to the system (\ref{HARPSysInf}) as
  \begin{gather*}
    F(X) = \lambda(X)\dpr{P(X)}{X}, \\
    Q(P(X)) = \frac1{\lambda(X)}.
  \end{gather*}
\end{Coroll}
See appendix~1 for the proofs of these corollaries.

Now we move to the concept of separability. Assume that we have a
group of goods split into $K+1$ subgroups and there are $m_k$ goods
in $k$-th subgroup (still $m$ is the total number of goods so that
$m = \sum_{k=1}^{K+1}m_k$). Let the elements of the consumption
vector $X$ be rearranged to reflect this split. Namely, assume that
$X = (X_1',\ldots,X_K', X_{K+1}')'$, where $X_k$ is the consumption
vector for the goods from $k$-th subgroup. Let the inverse demand
functions $P(X)$ be rearranged in the same manner:
\begin{gather*}
  P(X) = \left((P_1(X))',\ldots,(P_K(X))',
  (P_{K+1}(X))'\right)'.
\end{gather*}
The subgroups $1,\ldots,K$ are said to be separable from the
original group if the group is rationalizable and the rationalizing
function is represented as
\begin{gather*}
  F(X) = F_0(F_1(X_1),\ldots,F_K(X_K),X_{K+1}),
\end{gather*}
where all the functions $F_0$, $F_1$, $\ldots$, $F_K$ are from the
same functional class $\Phi$.

\setcounter{TheoremPrime}{\value{Theorem}}

\begin{Theorem}{\citep{VratenkovShananin1991Eng}}\label{HARPSep}
  Let the inverse demand function $P(X)$ be rationalizable in $\Phi_H$
  and the subgroups $1,\ldots,K$ be separable from the original
  group so that the rationalizing function is represented as
  \begin{gather*}
    F(X) = F_0(F_1(X_1),\ldots,F_K(X_K),X_{K+1}),
  \end{gather*}
  where $F_0$, $F_1$, $\ldots$, $F_K \in \Phi_H$. Then the Yang transform of $F(X)$
  \begin{gather*}
    Q(P) = \inf\left\{\left.\frac{\dpr{P}{X}}{F(X)}\right|X\geqslant
    0, F(X) > 0\right\}
  \end{gather*}
  may be represented as
  \begin{gather*}
    Q(P) = Q_0(Q_1(P_1),\ldots,Q_K(P_K),P_{K+1}),
  \end{gather*}
  where $P_k$ is price vector for the goods from the $k$-th subgroup,
  $Q_0(q_1, q_2, \ldots, q_K$, $P_{K+1})$ is the optimal value of
  objective function in the following optimization problem
  \begin{align*}
    & \frac{\sum\limits_{k=1}^Tq_kY_k +
    \dpr{P_{K+1}}{X_{K+1}}}{F_0(Y_1, Y_2, \ldots, Y_K, X_{K+1})} \rightarrow
    \min, \\
    & \text{\textup{s.t.} } (Y_1, Y_2, \ldots, Y_K, X_{K+1}) \in R_+^{K+m_{K+1}}, \\
    & \phantom{\text{\textup{s.t.} }} F_0(Y_1, Y_2, \ldots, Y_K, X_{K+1}) >
    0,
  \end{align*}
  and
  \begin{gather*}
    Q_k(P_k) = \inf\left\{\left. \frac{\dpr{P_k}{X_k}}{F_k(X_k)} \right|
    X_k \in \R_+^{m_k}, F_k(X_k) > 0 \right\} (k = 1, 2, \ldots, K).
  \end{gather*}
\end{Theorem}

\setcounter{CorollPrime}{\value{Coroll}}

\begin{Coroll} \label{VSCoroll}
  Let $P_1(X), \ldots, P_K(X)$ be rationalizable in $\Phi_H$. Let
  $Q_k(P_k)$ and $F_k(X_k)$ $(k=1,\ldots,K)$ be Kon\"us-Divisia price and consumption indices that
  is for all $k=1\ldots,K$
  \begin{align}
    & Q_k(P_k)F_k(X_k) \leqslant \dpr{P_k}{X_k} & & \forall P_k\in\R^{m_k}_+,
    X_k\in\R_+^{m_k}, \label{VSCoroll1} \\
    & Q_k(P_k(X))F_k(X_k) = \dpr{P_k(X)}{X_k} & & \forall X_k
    \in \R_+^{m_k}. \label{VSCoroll2}
  \end{align}
  Let the inverse
  demand functions
  \begin{gather*}
    \mathcal{P}(X) = (Q_1(P_1(X)),\ldots,Q_K(P_K(X)),P_{K+1}(X))
  \end{gather*}
  be rationalizable in $\Phi_H$ that is there exist functions
  $F_0(y_1,\ldots,y_K,X_{K+1})$ and $Q_0(q_1$, $\ldots$, $q_K$, $P_{K+1})$ from $\Phi_H$ such that
  \begin{gather}
    Q_0(q_1,\ldots,q_K,P_{K+1})F_0(y_1,\ldots,y_K,X_{K+1}) \leqslant \sum_{k=1}^Kq_ky_k +
    \dpr{P_{K+1}}{X_{K+1}}, \label{VSCoroll3}
  \end{gather}
  for all $(q_1,\ldots,q_K) \in \R_+^K$, $(y_1,\ldots,y_K) \in \R_+^K$,
  $P_{K+1}\in\R^{m_{K+1}}_+$, $X_{K+1}\in\R_+^{m_{K+1}}$, and
  \begin{multline}
    Q_0(Q_1(P_1(X)),\ldots,Q_K(P_K(X)),P_{K+1}(X))
    \\
    \times F_0(F_1(X_1),\ldots,F(X_K),X_{K+1}) =
    \sum_{k=1}^{K+1}\dpr{P_k(X)}{X_k}, \label{VSCoroll4}
  \end{multline}
  for all $X \in \R_+^{m}$. Then the inverse demand functions $P(X)$ are rationalizable in $\Phi_H$. Moreover, the
  Kon\"us-Divisia price and consumption indices for the whole group of
  goods are represented as
  \begin{gather}
    F(X) = F_0(F_1(X_1), \ldots, F_K(X_K), X_{K+1}), \label{VSCoroll5} \\
    Q(P(X)) = Q_0(Q_1(P_1(X)), \ldots, Q_K(P_K(X)), P_{K+1}(X)). \label{VSCoroll6}
  \end{gather}
\end{Coroll}
See appendix~1 for the proof.

These results may easily be reformulated for trade statistics rather
than inverse demand functions. Namely, let $TS = \{(P^t,
X^t)\}_{t=1}^T$ be trade statistics for the whole group of goods
with $X^t = ((X^t_1)', \ldots (X_K^t)', (X_{K+1}^t)')'$ and $P^t =
((P^t_1)', \ldots, (P_K^t)', (P_{K+1}^t)')'$. Assume that trade
statistics for the first $K$ subgroups satisfy HARP. Denote the
corresponding Kon\"us-Divisia price and consumption indices as
$Q_k(P^t_k)$ and $F_k(X_k^t)$. Construct new trade statistics
$\mathcal{TS} = \{(\mathcal{P}^t, \mathcal{X}^t)\}_{t=1}^T$ where
\begin{gather*}
  \mathcal{P}^t = (Q_1(P_1^t),\ldots, Q_K(P_K^t), P^t_{K+1}), \quad
  \mathcal{X}^t = (F_1(X_1^t),\ldots, F_K(X_K^t), X^t_{K+1}).
\end{gather*}

\begin{TheoremPrime}{\citep{VratenkovShananin1991Eng}}
  If the trade statistics $TS$ is rationalizable in $\Phi_H$ and
  corresponding Kon\"us-Divisia consumption indices are represented as
  \begin{gather*}
    F(X^t) = F_0(F_1(X_1^t),\ldots,F_K(X_K^t),X_{K+1}^t).
  \end{gather*}
  Then the corresponding Kon\"us-Divisia price indices may be
  represented as
  \begin{gather*}
    Q(P^t) = Q_0(Q_1(P_1^t),\ldots,Q_K(P_K^t),P_{K+1}^t),
  \end{gather*}
  where $Q_0(q_1^t, q_2^t, \ldots, q_K^t$, $P_{K+1}^t)$ is the optimal value of
  objective function in the following optimization problem
  \begin{align*}
    & \frac{\sum\limits_{k=1}^Tq_k^tY_k +
    \dpr{P^t_{K+1}}{X_{K+1}}}{F_0(Y_1, Y_2, \ldots, Y_K, X_{K+1})} \rightarrow
    \min, \\
    & \text{\textup{s.t.} } (Y_1, Y_2, \ldots, Y_K, X_{K+1}) \in R_+^{K+m_{K+1}}, \\
    & \phantom{\text{\textup{s.t.} }} F_0(Y_1, Y_2, \ldots, Y_K, X_{K+1}) >
    0,
  \end{align*}
  and
  \begin{gather*}
    Q_k(P_k^t) = \inf\left\{\left. \frac{\dpr{P_k^t}{X_k}}{F_k(X_k)} \right|
    X_k \in \R_+^{m_k}, F_k(X_k) > 0 \right\} (k = 1, 2, \ldots, K).
  \end{gather*}
\end{TheoremPrime}

\begin{CorollPrime}{\citep{VratenkovShananin1991Eng}}\label{CorVS}
  If the trade statistics $\mathcal{TS}$ satisfies HARP then the
  original trade statistics $TS$ also satisfies HARP. Moreover,
  Kon\"us-Divisia index numbers $F(X^t)$ and $Q(P^t)$ may be
  represented as
  \begin{gather*}
    F(X^t) = F_0(F_1(X_1^t), \ldots, F_K(X_K^t), X_{K+1}^t), \\
    Q(P^t) = Q_0(Q_1(P_1^t), \ldots, Q_K(P_K^t), P_{K+1}^t).
  \end{gather*}
\end{CorollPrime}

The corollary \ref{CorVS} shows that the Kon\"us-Divisia indices of
the original group are constructed from the Kon\"us-Divisia indices
of the subgroups. The separable subgroups may be further split into
smaller subgroups. When Kon\"us-Divisia index numbers are
constructed in this way we say that we have a hierarchy of
Kon\"us-Divisia index numbers. It may be represented as a graph in
which the vertices correspond to groups of goods and edges
correspond to the relation of separability. Such graph is one way of
visualizing the structure of consumer demand. An example of such
visualization is shown in subsection~\ref{HungTS}. In similar way
one may build a hierarchy of economic indices based on
GARP\footnote{In this case we mean only consumption indices as there
is no analog of theorem \ref{HARPSep} for rationalizability in
$\Phi_G$ showing the duality between consumption and price
indices.}.

\subsection{Irrationality indices}

The trade statistics may fail to satisfy HARP or GARP because of
errors in measurement of consumption and price data. This
observation suggests that one needs some quantitative measure of the
degree to what the trade statistics fails HARP or GARP.

Here we describe one approach to such quantitative measure. First,
let us introduce parametric families of axioms of revealed
preference. We say that the trade statistics satisfies
GARP($\omega$) if for all $t$ and $s$ the existence of
$t_1,\ldots,t_k$ such that $\dpr{P^t}{X^t}\geqslant \omega
\dpr{P^t}{X^{t_1}}$, $\dpr{P^{t_1}}{X^{t_1}}\geqslant \omega
\dpr{P^{t_1}}{X^{t_2}}$, $\ldots$, $\dpr{P^{t_k}}{X^{t_k}}\geqslant
\omega \dpr{P^{t_k}}{X^s}$ implies the inequality $\dpr{P^s}{X^s}
\leqslant \omega \dpr{P^s}{X^t}$. We say that the trade statistics
satisfies HARP($\omega$) if for all subsets of indices
$\{t_1,\ldots,t_k\}$ from $\{1,\ldots,T\}$ with $t_i \ne t_{i+1}$
for all $i \in \{1,\ldots, k-1\}$ the following inequality is
satisfied:
\begin{gather*}
  \dpr{P^{t_1}}{X^{t_2}}\dpr{P^{t_2}}{X^{t_3}}\ldots
  \dpr{P^{t_k}}{X^{t_1}} \geqslant \frac1{\omega^k}\dpr{P^{t_1}}{X^{t_1}}\dpr{P^{t_2}}{X^{t_2}}\ldots\dpr{P^{t_k}}{X^{t_k}}.
\end{gather*}

The additional requirement that $t_i \ne t_{i+1}$ is added in order
to make it possible for the trade statistics to satisfy
HARP($\omega$) with $\omega < 1$. We did not put this requirement in
the previous subsection because it would not have affected whether
certain trade statistics satisfied HARP or not. The Afriat-Varian
theorems from the previous subsection may be generalized in the
following way.

\begin{Theorem}\label{GGARPTh}
  The following statements are equivalent
  \begin{itemize}
    \item[1)] trade statistics $\{(P^t, X^t)\}_{t=1}^T$ satisfies GARP($\omega$).
    \item[2)] There exist numbers $U^t$, $\lambda^t > 0$,
    $(t=1,\ldots,T)$, such that
    \begin{gather}
      U^t \leqslant U^s + \lambda^s\left(\omega \dpr{P^s}{X^t} - \dpr{P^s}{X^s}\right), \quad \forall
      t,s = 1,\ldots,T, \quad t \ne s. \label{GGARPSys}
    \end{gather}
  \end{itemize}
\end{Theorem}

The proof is given in appendix~1.

\begin{Theorem}[\cite{Houtman1995}]
  The following statements are equivalent
  \begin{itemize}
    \item[1)] trade statistics $\{(P^t, X^t)\}_{t=1}^T$ satisfies HARP($\omega$);
    \item[2)] there exist numbers $\lambda^t > 0$ $(t =
    1,\ldots,T)$, such that
    \begin{gather}
      \omega\lambda^t\dpr{P^t}{X^s} \geqslant \lambda^s\dpr{P^s}{X^s}, \quad \forall t,s =
      1,\ldots,T, \quad t \ne s; \label{GHARPSys}
    \end{gather}
  \end{itemize}
\end{Theorem}

One may construct Kon\"us-Divisia index numbers using a solution to
the system (\ref{GHARPSys}) by formulas from the previous
subsection.

Now we may introduce the quantitative measures of
irrationalizability of trade statistics. For rationalizability in
$\Phi_G$ the measure is given by GARP irrationality index
\begin{gather*}
  \omega_G = \min\{\omega \mid \text{trade statistics satisfies
  GARP($\omega$)} \}.
\end{gather*}

This approach is close to that of \citep{Afriat1973} where a
critical efficiency level\footnote{It is equal to
$\frac1{\omega_G}.$} was used to measure the degree of inconsistency
of trade statistics with GARP. For rationalizability in $\Phi_H$ the
measure is given by HARP irrationality index
\begin{gather*}
  \omega_H = \min\{\omega \mid \text{trade statistics satisfies
  HARP($\omega$)} \}.
\end{gather*}

The method of constructing Kon\"us-Divisia index numbers using a
solution of the system (\ref{GHARPSys}) with $\omega = \omega_H$ is
called generalized nonparametric one (see
\citep{PospelovaShananin1998Eng}).

The theorem \ref{HARPandGARP} also may be generalized to the case of
$\omega \ne 1$.
\begin{Theorem} \label{GHARPandGGARP}
  The following statements are equivalent:
  \begin{itemize}
    \item[1)] trade statistics $\{(P^t, X^t)\}_{t=1}^T$ satisfies HARP($\omega$);
    \item[2)] trade statistics $\{(P^t, \mu^tX^t)\}_{t=1}^T$ satisfies
    GARP($\omega$) for any positive values of~$\mu^t$ $(t \in \{1,\ldots,T\})$.
  \end{itemize}
\end{Theorem}
See appendix~1 for the proof.

%There were other approaches to measure the degree of consistency of
%the trade statistics with GARP. One of them, suggested in
%\citep{Varian1990}, is to introduce a set of efficiency indices
%$\{\omega^t\}_{t=1}^T$, one for each time period, and then find the
%vector of irrationality indices as the solution to the quadratic
%programming problem
%\begin{align*}
%  & \sum_{t=1}^T(\omega^t - 1)^2 \rightarrow \min \\
%\end{align*}
%subject to the restriction that the trade statistics
%$\{(P^t,X^t)\}_{t=1}^T$ satisfies GARP($\omega^1$, $\ldots$,
%$\omega^T$), that is for any $t$ and $s$ the existence of $t_1$,
%$\ldots$, $t_k$ such that $P^tX^t\geqslant \omega^tP^tX^{t_1}$,
%$P^{t_1}X^{t_1}\geqslant \omega^{t_1}P^{t_1}X^{t_2}$, $\ldots$,
%$P^{t_k}X^{t_k} \geqslant \omega^{t_k}P^{t_k}X^s$ implies $P^sX^s <
%\omega^s P^sX^t$. Other approaches are reviewed in
%\citep{Gross1995}.

%For HARP we suggest a different approach close to one suggested in
%\citep{Varian1990}. Namely, we suggest to compute the vector of
%irrationality indices as the part of the solution to the following
%optimization problem:
%\begin{align*}
%  & \sum_{t=1}^T\log(\omega^t)
%  \rightarrow \min\limits_{\{\omega^t, \lambda^t\}_{t=1}^T} \\
%  & \text{s.t.} \log(\omega^t) + \log(\lambda^t) - \log(\lambda^s)
%  \geqslant \log\left(\frac{P^sX^s}{P^tX^s}\right) \quad (t,s \in
%  \{1,\ldots,T\}, t\ne s).
%\end{align*}
%
%This problem is simpler than one suggested in \citep{Varian1990}
%because it is a linear programming problem and it may be effectively
%solved by simplex method.

\subsection{Forecasting}\label{SSForecasting}

Assume that we have some trade statistics $\{(P^t, X^t)\}_{t=1}^T$.
There are two forecasting problems which may be solved with revealed
preference theory. The first one is to predict $P^{T+1}$
given~$X^{T+1}$. The second one is to predict $X^{T+1}$
given~$P^{T+1}$. Since both problems have similar solution, we
describe the details only for the second problem.

Assume that the trade statistics satisfies GARP($\omega$). Then we
may construct the set of forecasted demand vectors as
\begin{gather*}
  K_G^{\omega}(P^{T+1}) = \left\{X^{T+1} \left| \{(P^t,X^t)\}_{t=1}^{T+1} \text{ satisfies
  GARP($\omega$)}\right\}\right. .
\end{gather*}
This approach is a generalization of that for GARP(1) suggested in
\citep{Varian1982a}.

%\begin{Theorem}
%  Let $N(P^{T+1}) = \{\tau\mid\exists \{t_1,\ldots,t_k\}\subset
%  \{1,\ldots,T\} : P^{\tau}X^{\tau}\geqslant \omega
%  P^{t_1}X^{\tau}$,
%  $P^{t_1}X^{t_1}\geqslant \omega P^{t_2}X^{t_1}, \ldots,
%  P^{t_k}X^{t_k}\geqslant \omega P^{T+1}X^{t_k} \}$. Then,
%  \begin{gather*}
%    K_G^{\omega}(P^{T+1}) = \{ X \mid P^{T+1}X \leqslant \omega P^{\tau}X, \forall
%    \tau \in N(P^{T+1})\}.
%  \end{gather*}
%\end{Theorem}

%This theorem is a corollary of theorem \ref{GGARPTh}. In order to
%find $N(P^{T+1})$ one one starts with construction of transitive
%closure $\mathit{R}^*_P$ of direct revealed preference relation
%defined on $\{P^t\}_{t=1}^T \times \{P^t\}_{t=1}^T$ by
%\begin{gather*}
%  P^t\mathit{R}_P(\omega) P^s \Leftrightarrow P^tX^t \geqslant
%  \omega P^sX^t.
%\end{gather*}
%Then one finds $n(P^{T+1}) = \{\tau \mid \tau \in \{1,\ldots,T\},
%P^{\tau}X^{\tau} \geqslant P^{\tau}X^{T+1}\}$ and
%\begin{gather*}
%  N(P^{T+1}) = \mathop{\bigcup}\limits_{t\in n(P^{T+1})}\{\tau \mid
%  P^{\tau}\mathit{R}^*_P P^t\}.
%\end{gather*}

In the similar way we may construct the set of forecasted demand
vectors assuming that the trade statistics satisfies HARP($\omega$)
as
\begin{gather*}
  K_H^{\omega}(P^{T+1}) = \left\{X^{T+1} \left| \{(P^t,X^t)\}_{t=1}^{T+1} \text{ satisfies
  HARP($\omega$)}\right\}\right. .
\end{gather*}

\begin{Theorem}{\citep{GrebShan08Eng}}

  Assume that the trade statistics $\{(P^t, X^t)\}_{t=1}^T$ satisfies HARP($\omega$) with $\omega \geqslant 1$. Let
  \begin{gather*}
    C_{ts}^*(\omega) =
    \max\left\{\omega^{-k-1}C_{tt_1}C_{t_1t_2}\ldots C_{t_{k-1}t_k}C_{t_ks} \mid \{t_1,\ldots,t_k\} \subset
    \{1,\ldots,T\},\,k\geqslant0\right\},
  \end{gather*}
  and
  \begin{gather*}
    \gamma_s\left(P^{T+1}, \omega\right) =
    \min\limits_{t\in\{1,\ldots,T\}}\left\{\frac{\omega^2}{C^*_{ts}(\omega)}
    \frac{\dpr{P^{T+1}}{X^{t}}}{{P^t}{X^t}}\right\}.
  \end{gather*}
  Then
  \begin{gather*}
    K_H^{\omega}(P^{T+1}) = \left\{X \in \mathbb{R}^m_+ \mid \gamma_{s}\left(P^{T+1},\omega\right)
  \dpr{P^{s}}{X} \geqslant \dpr{P^{T+1}}{X} \, \forall s \in
  \{1,\ldots,T\}\right\}.
  \end{gather*}
\end{Theorem}
The values $C^*_{ts}(\omega)$ may be effectively computed in
$O(T^3)$ operations by means of Floyd-Warshall algorithm.

We may fix not only new price vector at the value $P^{T+1}$ but also
the expenditure at some level $x_{T+1}$. Then, the corresponding
forecasting sets are given by
\begin{gather*}
  K_G^{\omega}(P^{T+1}, x_{T+1}) = K_G^{\omega}(P^{T+1}) \cap \left\{X \in
  \mathbb{R}^m_+ \mid \dpr{P^{T+1}}{X} = x_{T+1}\right\}, \\
  K_H^{\omega}(P^{T+1}, x_{T+1}) = K_H^{\omega}(P^{T+1}) \cap \left\{X \in
  \mathbb{R}^m_+ \mid \dpr{P^{T+1}}{X} = x_{T+1}\right\}.
\end{gather*}
The analytical expression for $K_G^{1}(P^{T+1}, 1)$ may be found in
\citep{Varian1982a}.

One may propose forecasting sets based on the concept of
rationalizability of Engel curves. Let the set of Engel curves
$\{q^t(\cdot)\}_{t=1}^T$ corresponding to price vectors
$\{P^t\}_{t=1}^T$ satisfy GARP($\omega$). The latter means that for
any set of individual expenditures $\{x_t\}_{t=1}^T$ the trade
statistics $\left\{(P^t, q^t(x_t))\right\}_{t=1}^T$ satisfies
GARP($\omega$). Then one may suggest to build the forecasting set
for demand at a fixed price vector $P^{T+1}$ as
\begin{gather*}
  \widetilde{K}_G^{\omega}(P^{T+1}) = \left\{X \in \mathbb{R}^m_+ \left|
  \begin{array}{c}
    \{(P^t, q^t(x_t))\}_{t=1}^{T} \cup \{(P^{T+1},X)\} \text{ satisfies GARP($\omega$)} \\
    \text{for any $x_1 > 0,
  \ldots, x_T > 0$}
  \end{array}
  \right.\right\}.
\end{gather*}
Denote the set $\widetilde{K}_G^{\omega}(P^{T+1}) \cap \{X \in
\mathbb{R}^m_+ \mid \dpr{P^{T+1}}{X} = x_{T+1}\}$ as
$\widetilde{K}_G^{\omega}(P^{T+1}, x_{T+1})$.

Let us introduce also the set $\widetilde{K}_H^{\omega}(P^{T+1})$ as
\begin{gather*}
  \widetilde{K}_H^{\omega}(P^{T+1}) = \left\{X \in \mathbb{R}^m_+ \left|
  \begin{array}{c}
    \{(P^t, q^t(x_t))\}_{t=1}^{T} \cup \{(P^{T+1},X)\} \text{ satisfies HARP($\omega$)} \\
    \text{for any $x_1 > 0,
  \ldots, x_T > 0$}
  \end{array}
  \right.\right\}.
\end{gather*}
It is worth noting that $\widetilde{K}_H^{1}(P^{T+1}) =
K_H^1(P^{T+1})$. That is because if a trade statistics $\{(P^t,
q^t(x_t) )\}_{t=1}^T$ satisfies HARP for some positive total
expenditures $x_1,\ldots,x_K$ then it is rationalizable with
positively-homogeneous utility function. This implies that the Engle
curves $q^t(\cdot)$ are rays. Therefore, choosing total expenditures
different from $x_1,\ldots,x_K$ is equivalent to change of scales in
consumption data. Since HARP is invariant to such changes, the trade
statistics $\{(P^t, q^t(x_t'))\}_{t=1}^T$ satisfies HARP for any
positive total expenditures $x_1',\ldots,x_K'$.

A very close approach was studied by
\citep{BlundellBrowningCrawford2008}. The difference is that the
authors used Strong Axiom of Revealed Preference (SARP). It is
similar to GARP though is little bit more restrictive. The trade
statistics $\{(P^t, X^t)\}_{t=1}^T$ satisfies SARP if
$X^t\mathit{R}^* X^s$ and $X^t\ne X^s$ implies $\dpr{P^s}{X^s} <
\dpr{P^s}{X^t}$. The authors of \citep{BlundellBrowningCrawford2008}
consider the following support sets
\begin{multline*}
  S(P^{T+1}, x_{T+1}, \{x_t\}_{t=1}^T) = \\
  \left\{X^{T+1} \in \R^m_+ \left|
  \begin{array}{c}
    \dpr{P^{T+1}}{X^{T+1}} = x_{T+1}, \\
    \{(P^t, q^t(x_t))\}_{t=1}^{T} \cup \{(P^{T+1},X)\}
    \text{ satisfies SARP},
  \end{array}
  \right. \right\}.
\end{multline*}

Similarly, one may consider support sets based on GARP:
\begin{multline*}
  G(P^{T+1}, x_{T+1}, \{x_t\}_{t=1}^T) = \\
  \left\{X^{T+1} \in \R^m_+ \left|
  \begin{array}{c}
    \dpr{P^{T+1}}{X^{T+1}} = x_{T+1}, \\
    \{(P^t, q^t(x_t))\}_{t=1}^{T} \cup \{(P^{T+1},X)\}
    \text{ satisfies GARP},
  \end{array}
  \right. \right\}.
\end{multline*}

The authors define intersection demands $\Tilde X^t = q^t(\tilde
x_t)$ as
\begin{gather*}
  \dpr{P^{T+1}}{q^t(\tilde x_t)} = x_{T+1},
\end{gather*}
and claim that if demands are weakly normal\footnote{This means that
if $x > x'$ then $q^t(x) \geqslant q^t(x')$ for all $t =
\overline{1,T}$.} then for any $\{x_t\}_{t=1}^T$
\begin{gather}
  S(P^{T+1}, x_{T+1}, \{\tilde x_t\}_{t=1}^T) \subseteq S(P^{T+1}, x_{T+1},
  \{x_t\}_{t=1}^T). \label{BBCTheorem}
\end{gather}
Since the closure of $S(P^{T+1}, x_{T+1}, \{x_t\}_{t=1}^T)$ is equal
to $G(P^{T+1}, x_{T+1}, \{x_t\}_{t=1}^T)$, this implies that
$\widetilde{K}_G^1(P^{T+1}, x_{T+1}) = G(P^{T+1}, x_{T+1}, \{\tilde
x_t\}_{t=1}^T)$. This conclusion is incorrect. We provide
counterexample showing its incorrectness in Appendix~2.

The further decrease in the size of the set $K_H^{\omega}$ is
possible if we add the additional requirement for the trade
statistics to satisfy the Law of Demand. This law states that the
price index is decreasing function of the demand index. We say that
the trade statistics satisfies the Law of Demand if it may be
continued to the demand functions satisfying the Law of Demand.

Let us define $\widehat{K}_H^{\omega}$ by
\begin{gather*}
  \widehat{K}_H^{\omega} = \left\{X^{T+1} \left| \{(P^t,X^t)\}_{t=1}^{T+1} \text{ satisfies
  HARP($\omega$) and the Law of Demand}\right\} \right. .
\end{gather*}

The analytical expression for $\widehat{K}_H^{\omega}$ is not known
yet. However, \citep{GrebShan08Eng} provides external estimate for
this set.

\begin{Theorem}{\citep{GrebShan08Eng}}

  Assume that the trade statistics $\{(P^t,X^t)\}_{t=1}^T$ satisfies HARP($\omega$) with $\omega \geqslant 1$. Let
  \begin{gather}
    D_{st} = \max\left\{\frac{\dpr{P^t}{X^t}}{\omega \dpr{P^s}{X^t}}, \Theta\left(
    \frac{\dpr{P^s}{X^s}}{\dpr{P^s}{X^t}} - \omega\right)\right\},
  \end{gather}
  where
  \begin{gather*}
    \Theta(x) = \left\{
    \begin{array}{lc}
      1, & x > 0, \\
      0, & x \leqslant 0,
    \end{array}
    \right.
  \end{gather*}
  and
  \begin{gather*}
    \Delta_{st} = \max\limits_{\{t_1,\ldots,t_k\} \subset
    \{1,\ldots,T\}}\left\{D_{st_1}D_{t_kt}
    \prod\limits_{i=2}^kD_{t_{i-1}t_i}\right\}.
  \end{gather*}
  Then
  \begin{multline*}
    \widehat{K}_H^{\omega} \subset \left\{
    \vphantom{X^{T+1} \in \R_+^m \left| \max\left[
    \frac{ \dpr{P^{T+1}}{X^{T+1}}}{\omega \dpr{P^s}{X^{T+1}}}, \Theta\left( \frac{
    \dpr{P^s}{X^s}}{\dpr{P^s}{X^{T+1}}} - 1\right)\right] \times \max\left[\frac{
    \dpr{P^t}{X^t}}{\omega \dpr{P^{T+1}}{X^t}}, \Theta\left( \frac{
    \dpr{P^{T+1}}{X^{T+1}}}{\omega \dpr{P^{T+1}}{X^t}} - 1\right)\right] \leqslant
    \frac1{\Delta_{st}} ,\quad \forall s,t = \overline{1,T}\right.}
    X^{T+1} \in \R_+^m \left| \max\left[
    \frac{ \dpr{P^{T+1}}{X^{T+1}}}{\omega \dpr{P^s}{X^{T+1}}}, \Theta\left( \frac{
    \dpr{P^s}{X^s}}{\dpr{P^s}{X^{T+1}}} - 1\right)\right]\right.\right. \\
    \left. \times \max\left[\frac{
    \dpr{P^t}{X^t}}{\omega \dpr{P^{T+1}}{X^t}}, \Theta\left( \frac{
    \dpr{P^{T+1}}{X^{T+1}}}{\omega \dpr{P^{T+1}}{X^t}} - 1\right)\right] \leqslant
    \frac1{\Delta_{st}} ,\quad \forall s,t = \overline{1,T}
    \vphantom{X^{T+1} \in \R_+^m \left| \max\left[
    \frac{ \dpr{P^{T+1}}{X^{T+1}}}{\omega \dpr{P^s}{X^{T+1}}}, \Theta\left( \frac{
    \dpr{P^s}{X^s}}{\dpr{P^s}{X^{T+1}}} - 1\right)\right] \times \max\left[\frac{
    \dpr{P^t}{X^t}}{\omega \dpr{P^{T+1}}{X^t}}, \Theta\left( \frac{
    \dpr{P^{T+1}}{
    X^{T+1}}}{\omega \dpr{P^{T+1}}{X^t}} - 1\right)\right] \leqslant
    \frac1{\Delta_{st}} ,\quad \forall s,t = \overline{1,T}\right.}
    \right\}.
  \end{multline*}
\end{Theorem}

\section{Applications of Kon\"us-Divisia index numbers}\label{SApplics}

In this section we describe several empirical applications of
Kon\"us-Divisia index numbers introduced in \ref{HARPGARP}. We used
two data sets. The first one contains consumption and price data for
196 goods in Hungary for years 1975-1984 \citep{HungTS}. During
these years Hungarian economy was shifting from plan to market one
and it is interesting to analyze the changes in the structure of
consumer demand happened in this period. The second one contains
similar data for 106 goods in the Netherlands for years 1951-1977
\citep{NethTS}. These are post-war years. The both data sets covers
the period when middle class was emerging.

\subsection{Trade statistics of Hungary}\label{HungTS}

In figure \ref{MeatGraphs} we show prices and consumption for six
meat products from Hungarian trade statistics (in thin lines) as
well as Kon\"us-Divisia consumption and price indices (in thick
lines). We see that Kon\"us-Divisia index numbers are less volatile
compared to raw prices and consumption. This implies that they may
be forecasted more accurately than prices and consumption of goods
themselves.

\begin{figure}[!h]
  \begin{center}
    \includegraphics[width = \linewidth]{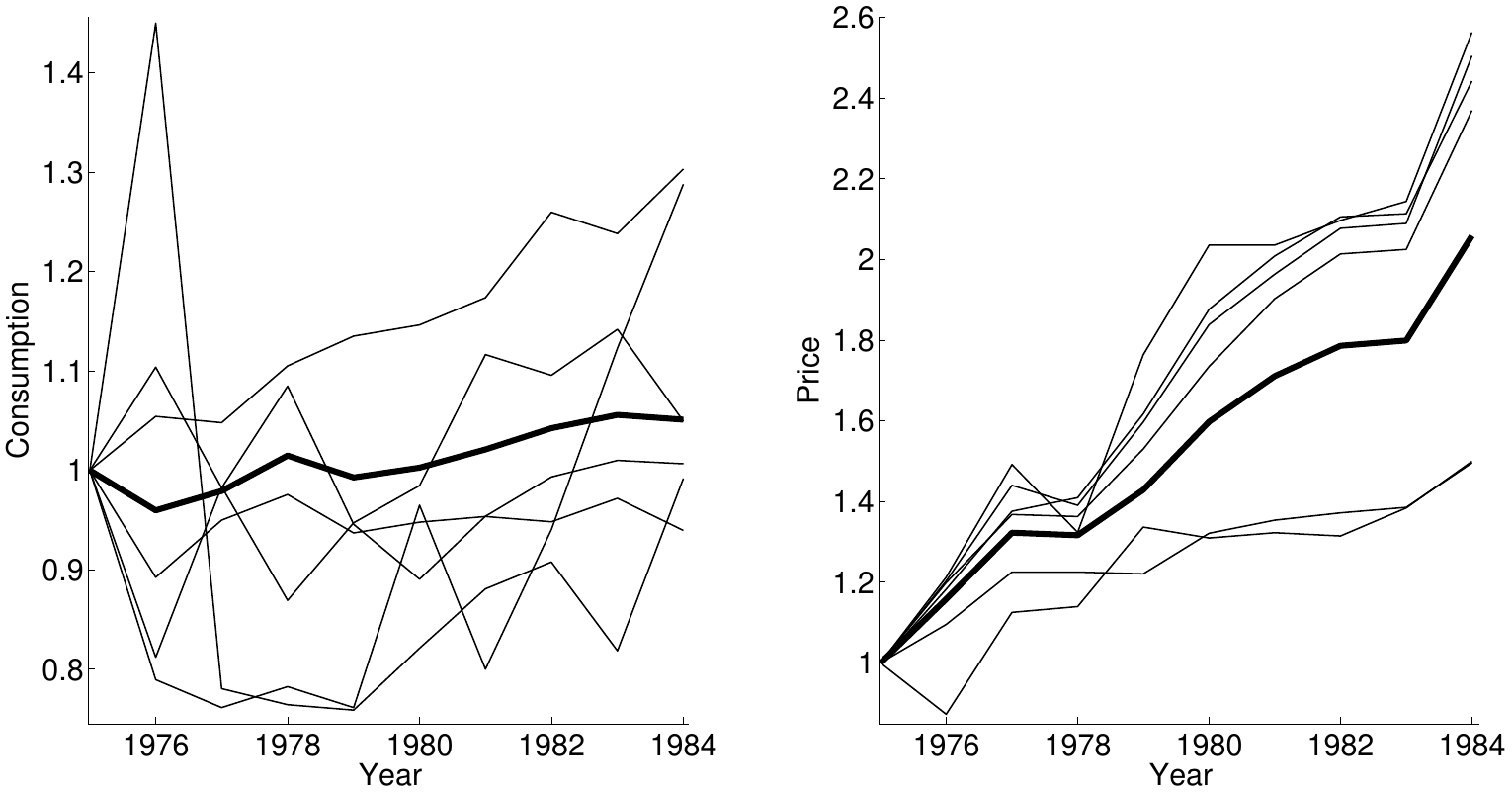}
    \caption{Prices, consumption, and Kon\"us-Divisia indices for group of meat product.}
    \label{MeatGraphs}
  \end{center}
\end{figure}

The commodities in Hungarian trade statistics are split by Hungary
trade experts into 11 groups shown in table \ref{HungGroups}. The
commodity groups differ by the duration of commodity usage. The
first three groups represent the everyday goods which are consumed
in 1 to 3 months. The goods from ''Clothes'' groups are used during
a year on average. All the rest groups contain goods and services of
durable consumption with characteristic duration of service of 5-10
years.

\begin{table}[!h]
  \begin{center}
  \caption{Commodity groups in the trade statistics of Hungary.}
  \label{HungGroups}
    \begin{tabular}{|c|l|c|}
      \hline
      Gr. number & Group name & Num. of goods \\
      \hline
      1 & Consumption goods & 49 \\
      2 & Beverages & 15 \\
      3 & Tobacco goods & 3 \\
      4 & Clothes & 31 \\
      5 & House service & 5 \\
      6 & Heating, home energy & 12 \\
      7 & Home equipment & 30 \\
      8 & Health care, hygiene & 7 \\
      9 & Transport, information & 11 \\
      10 & Education, culture, sport, recreation & 23 \\
      11 & Other goods & 10 \\
      \hline
    \end{tabular}
  \end{center}
\end{table}

The time series of Kon\"us-Divisia indices for the groups of all
goods, of everyday goods, and of goods and services of durable
consumption are shown in figure \ref{HungTreeGraphs}. We see that
the demand was shifting from everyday goods to goods and services of
durable consumption. The growth rate of the consumption of goods and
services of durable consumption was higher that of everyday goods.
This shift may be attributed to formation of middle class in Hungary
during the 1974-1984. This change in social structure lead to the
change of the hierarchy of economic indices.

\begin{figure}[!h]
  \begin{center}
    \includegraphics[width = \linewidth]{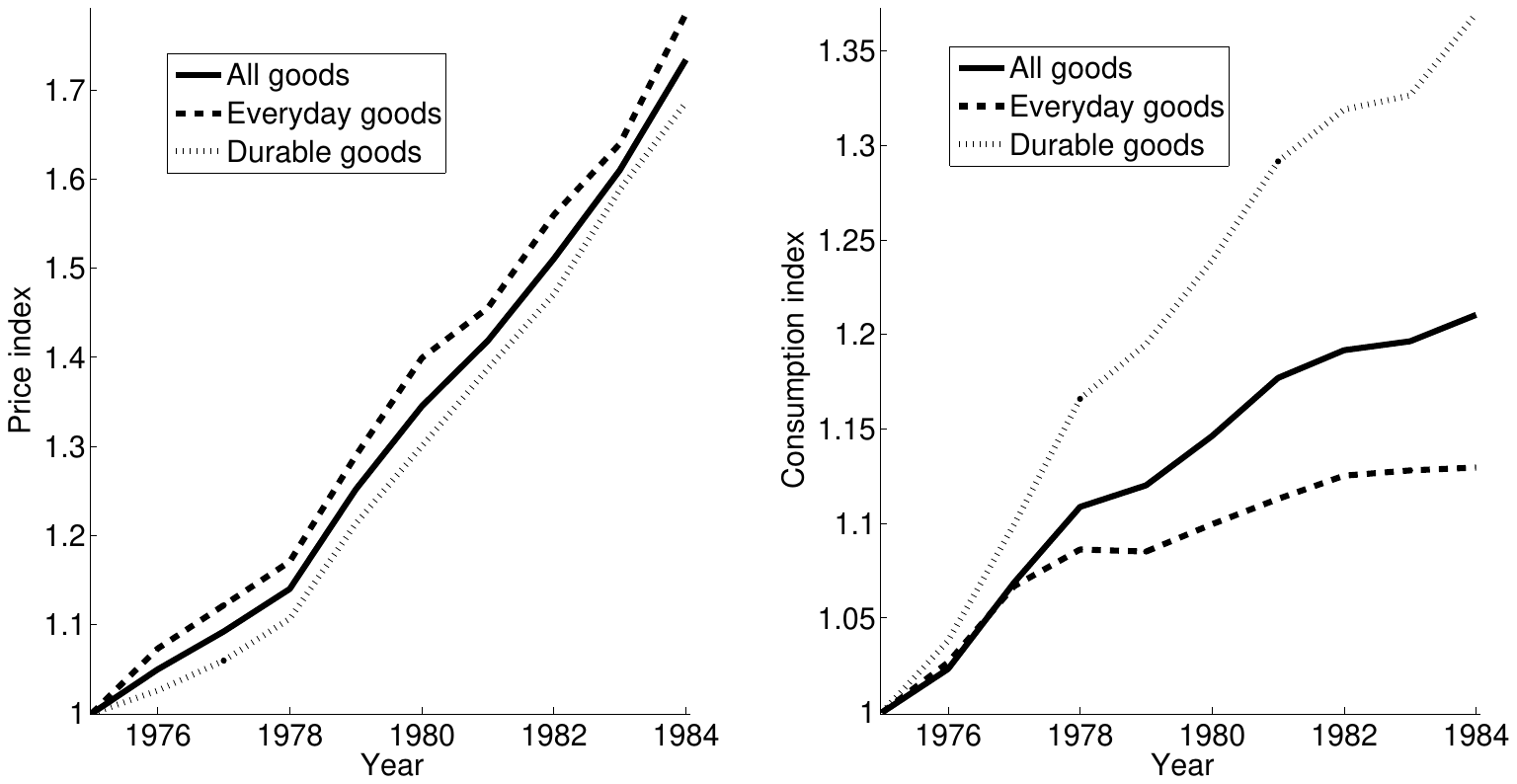}
    \caption{Segmentation in Hungarian trade statistics.}
    \label{HungTreeGraphs}
  \end{center}
\end{figure}

Out of the 11 groups selected by Hungary trade experts only the
first two satisfy HARP. The group formed by all goods also satisfies
HARP. The fact that the majority of the commodity groups selected by
Hungary trade experts does not satisfy HARP implies that such
classification does not reflect the structure of consumer demand in
Hungary during 1975-1984.

A different classification based on characteristic time of
consumption was suggested in \citep{VratenkovShananin1991Eng}. The
goods from the first three groups from table \ref{HungGroups} form a
group which satisfies HARP. These goods form new aggregated good
which is called as ''everyday goods'' in
\citep{VratenkovShananin1991Eng}. The goods from the groups from 5
to 11 in table \ref{HungGroups} also form a group satisfying HARP.
The corresponding aggregated good is named ''goods and services of
durable consumption'' in \citep{VratenkovShananin1991Eng}. The group
''Clothes'' does not satisfy HARP, but if did then we would have
another aggregated good called ''goods with medium time of
consumption''. However, the group with all the goods from the group
''Clothes'' and the aggregated good ''everyday goods'' does satisfy
HARP. The group formed by the group ''Clothes'' plus the aggregated
good ''goods and services of durable consumption'' also satisfies
HARP. Moreover, both these groups prove to be separable from the
group formed by all the goods. The latter group also satisfies HARP.
This structure of consumer demand was revealed in
\citep{VratenkovShananin1991Eng} and is shown in figure
\ref{HungTree}.

\begin{figure}[!h]
  \begin{center}
    \includegraphics[width = \linewidth]{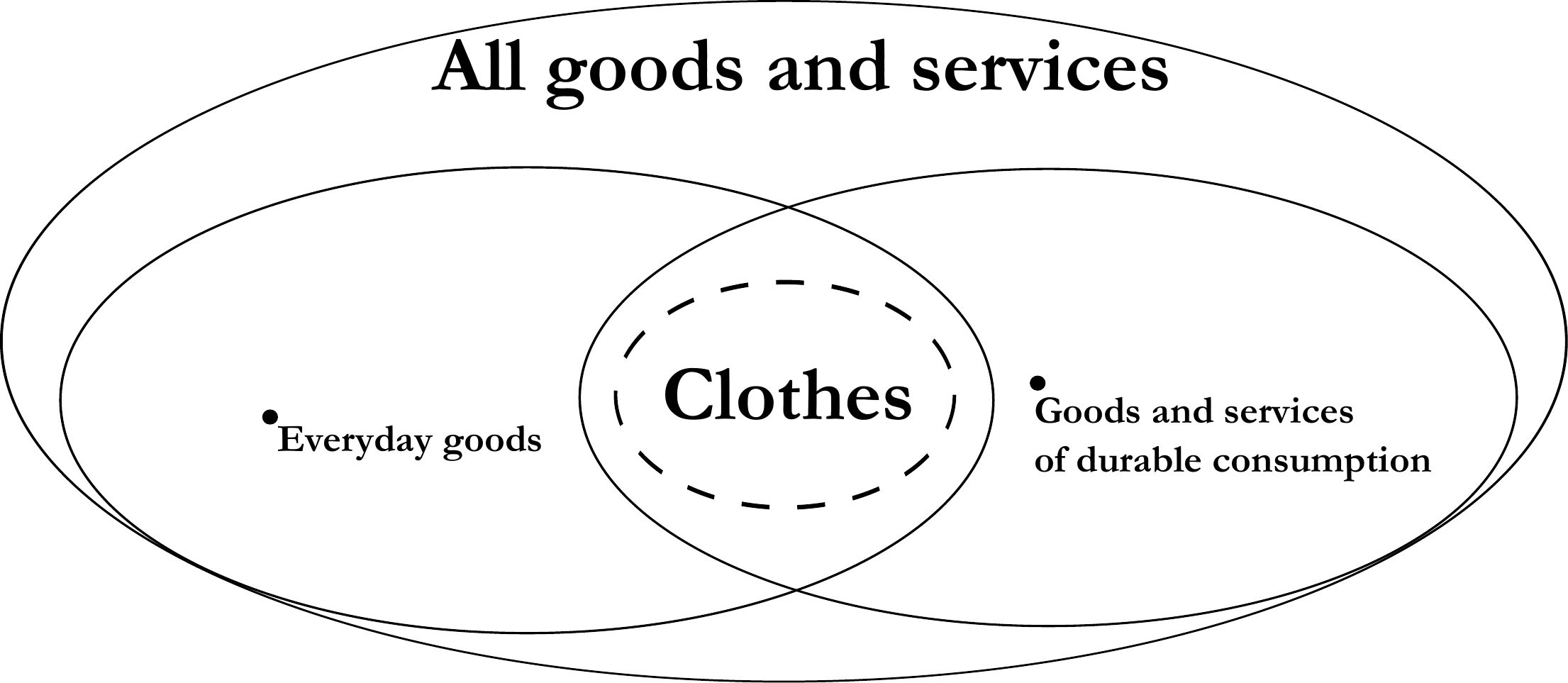}
    \caption{The structure of consumer demand in Hungary during 1975-1984.}
    \label{HungTree}
  \end{center}
\end{figure}

More details on studying this statistic may be found in
\citep{VratenkovShananin1991Eng, PospelovaShananin1998Eng,
KondrakovPospelovaShananin2010Eng}.

We used this statistics in order to asses the probability of having
a randomly chosen group satisfying GARP or HARP for various sizes of
random groups. More precisely, for each size of random groups in
range from 3 to 193 we generated 100000 random groups\footnote{The
number of different groups of sizes 2, 194, 195, and 196 is less
than 100000 so we considered all of them. We did not considered
groups consisting of a single good because any such group satisfies
both GARP and HARP.} of goods and computed the fraction of them
satisfying GARP and HARP. The values of computed fractions are
plotted against the sizes of random groups in figure
\ref{RandomGroups}. We see that for any size of random group the
probability that it satisfies GARP is close to 1. This means that
there will be a lot of groups which are complete in terms of
substitutes and compliments if we use GARP in order to build the
hierarchy of economic indices. This means that GARP allows almost
any structure of consumer demand and does not lead to revealing the
true one.

\begin{figure}[!h]
  \begin{center}
    \includegraphics[width = \linewidth]{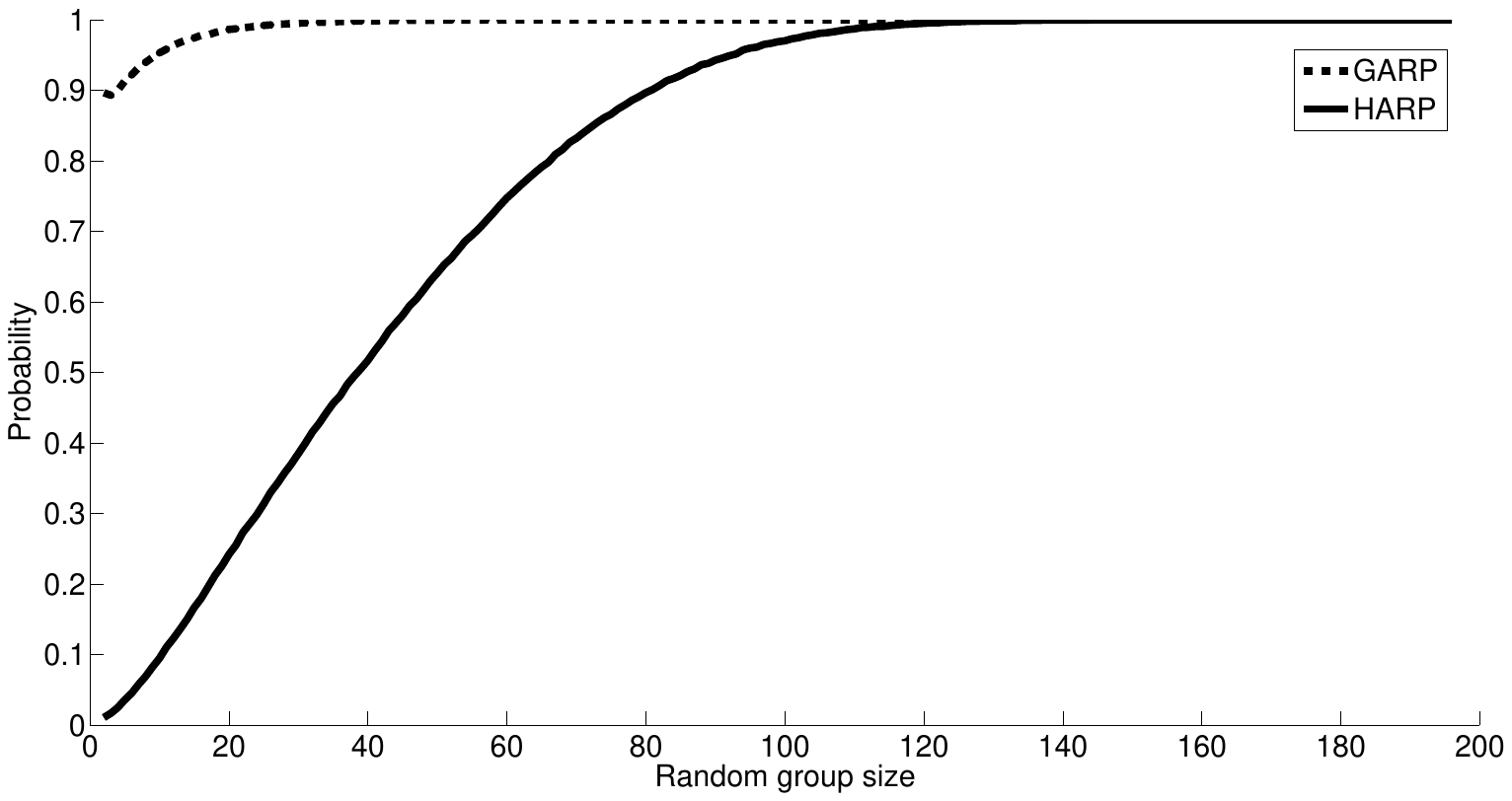}
    \caption{Probability of having random groups satisfying HARP and GARP on statistics of the Hungary.}
    \label{RandomGroups}
  \end{center}
\end{figure}

The situation is different for HARP. Although almost all large
groups of goods satisfy HARP\footnote{This is because the group
containing all the goods satisfies HARP and exclusion of small
number of goods leads only to slight violations of HARP in terms of
HARP irrationality index.}, the majority of small groups does not
satisfy HARP. This means that a group satisfying HARP is more
valuable and meaningful for the construction of the tree of economic
indices since there is more evidence in favor of it being an
independent market segment.

\subsection{Trade statistics of the Netherlands}

This trade statistics was analyzed in
\citep{PospelovaShananinBook1998Eng}. Dutch trade experts split all
the goods into five groups. None of them satisfy HARP. However the
group of all goods does satisfy HARP. The commodity groups with
corresponding irrationality indices are shown in table
\ref{NethGroups}.

\begin{table}[!h]
  \begin{center}
    \caption{Commodity groups in the trade statistics of the Netherlands.}
    \label{NethGroups}
    \begin{tabular}{|c|l|c|c|}
      \hline
      Gr. number & Goods & Num. of goods & $\omega_H$ \\
      \hline
      0 & All goods & 106 & 0.9998 \\
      \hline
      1 & Grocery, milk, meat, fish & 42 & 1.0002 \\
      \hline
      2 & Sweets, tobacco goods, beer & 12 & 1.0034 \\
      \hline
      3 & Clothes, cars, home equipment & 21 & 1.0073 \\
      \hline
      4 & Heating, water, medicine & 11 & 1.0002 \\
      \hline
      5 & Renting, service & 20 & 1.0167 \\
      \hline
    \end{tabular}
  \end{center}
\end{table}

Two ways of evaluating Kon\"us-Divisia indices for the group of all
goods were considered in \citep{PospelovaShananinBook1998Eng}. The
first way is to compute indices for the groups selected by trade
experts using generalized nonparametric method, then form new trade
statistics of the resulting Kon\"us-Divisia indices and evaluate
Kon\"us-Divisia indices for it. Denote price indices computed in
this way by $Q_1(P^t)$. The second way is to apply nonparametric
method directly to the group of all goods without intermediate
aggregation. Denote price indices computed in this way as
$Q_2(P^t)$.

The time series of these two price indices were obtained in
\citep{PospelovaShananinBook1998Eng} and are shown in figure
\ref{Neth}. We see that there is little difference between them.
Indeed, the maximum deviation between $Q_1(P^t)$ and $Q_2(P^t)$
is~1.76\%. This provides evidence in favor of robustness of
Kon\"us-Divisia indices to methods of aggregation.

\begin{figure}[!h]
  \begin{center}
    \includegraphics[width = \linewidth]{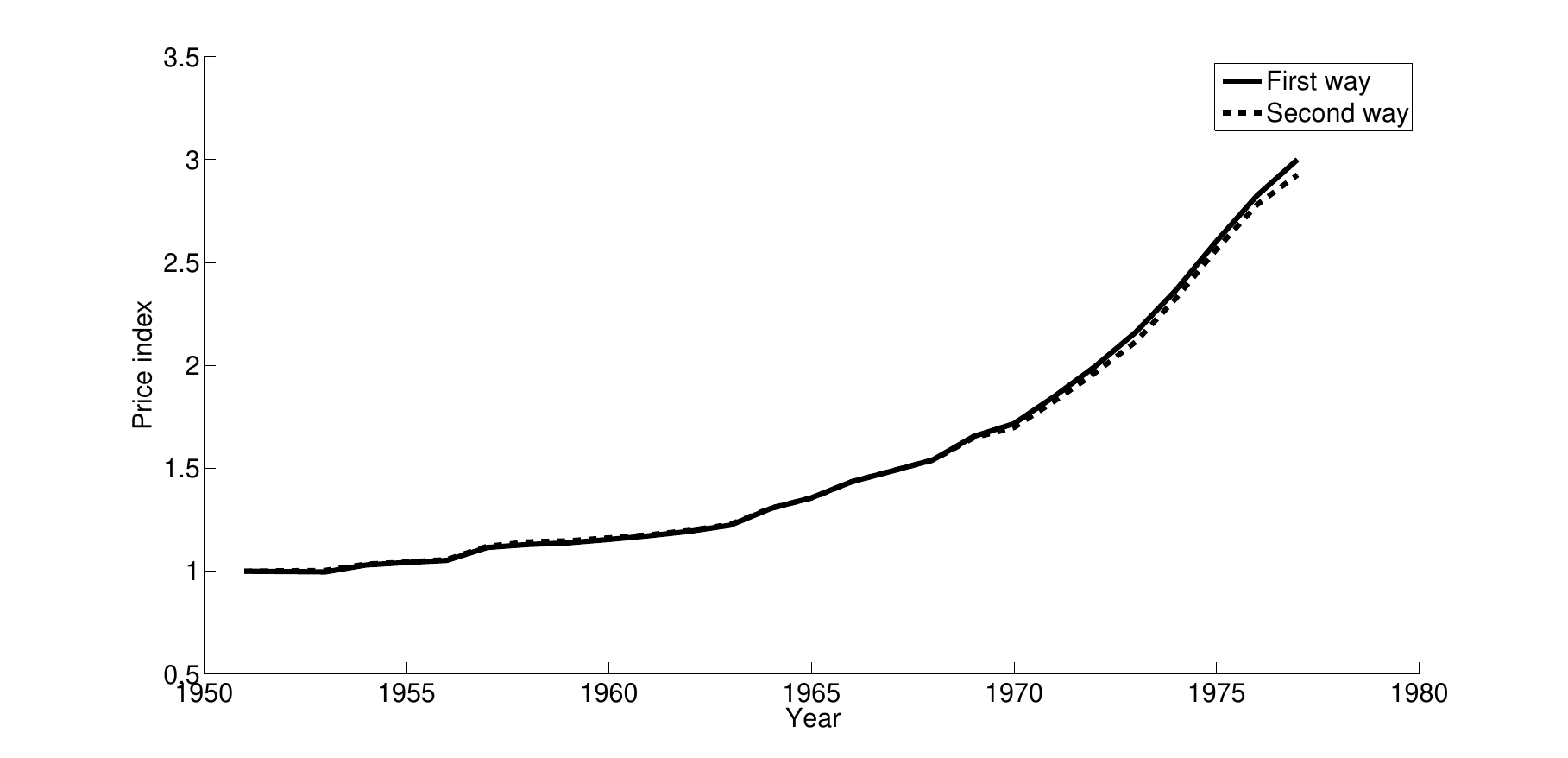}
    \caption{Kon\"us-Divisia price indices constructed in two different ways.}
    \label{Neth}
  \end{center}
\end{figure}

We repeated the empirical study of the probability of having random
group satisfying HARP or GARP from the previous subsection on the
statistics of the Netherlands. The estimated probabilities are shown
in figure \ref{RandomGroupsNeth}. Since in this statistics we have
more than twice time periods than in the statistics of Hungary we
expected the probabilities become lower. For HARP this is indeed so.
We see that the probability of having random group satisfying HARP
does not reach the value of one for any random group size below the
total number of goods. The probability of having a random group
satisfying GARP became lower for small random groups. However, it is
still extremely close to one for random groups of sizes greater than
30.

\begin{figure}[!h]
  \begin{center}
    \includegraphics[width = \linewidth]{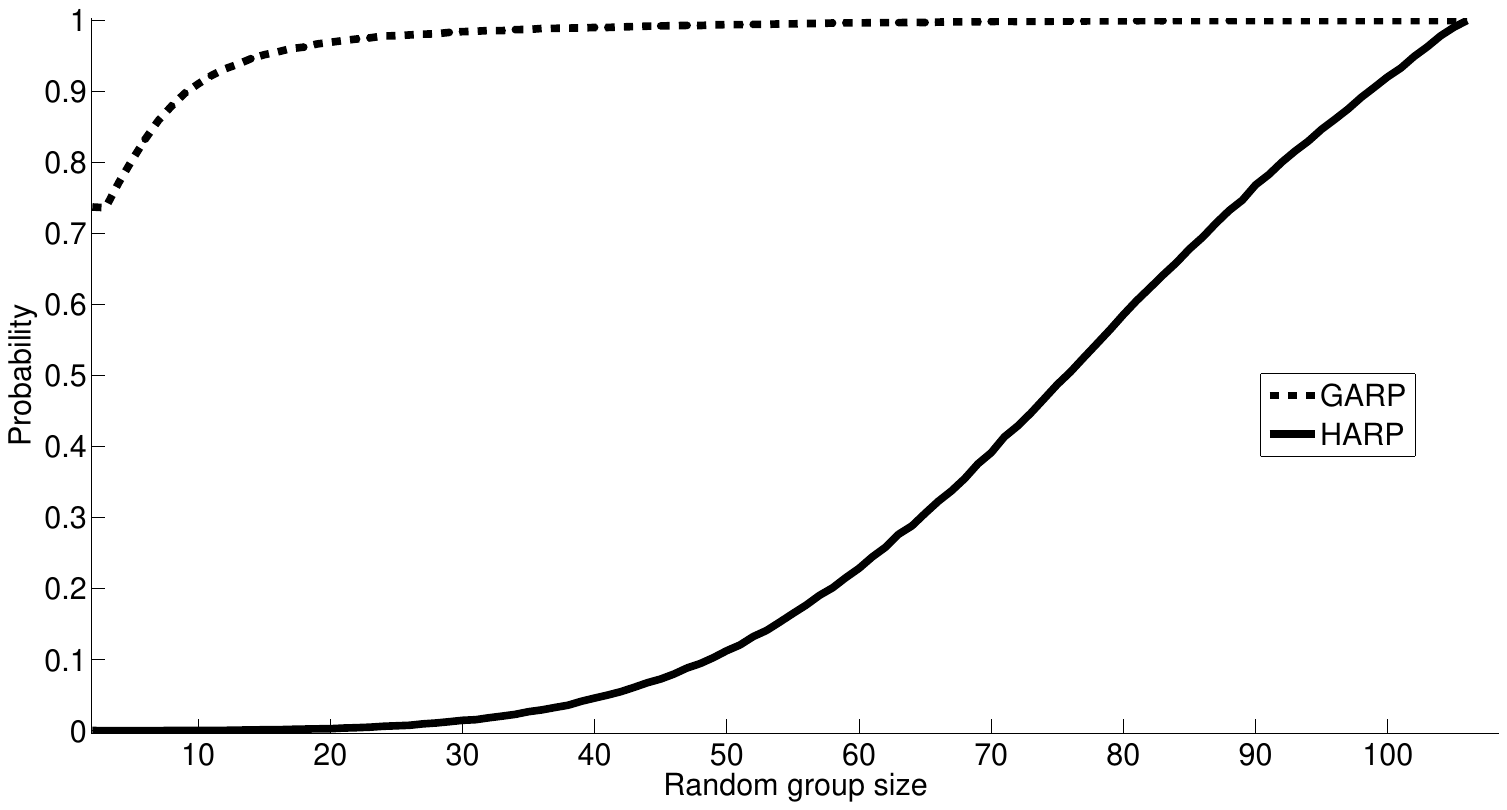}
    \caption{Probability of having random groups satisfying HARP and GARP on statistics of the Netherlands.}
    \label{RandomGroupsNeth}
  \end{center}
\end{figure}

\section{Power of HARP and GARP tests}\label{SPower}

In this section we assess the power of tests for consistency with
maximization behavior using the trade statistics of Hungary. The
general approach is the following. Consider a trade statistics
$\{(P^t, X^t)\}_{t=1}^T$. Denote random trade statistics with random
prices and original consumption as $\widetilde{TS} = \{(\Tilde P^t,
X^t)\}_{t=1}^T$. We generate $B$ realizations of $\widetilde{TS}$:
$\{\widetilde{TS}_b\}_{b=1}^B$. For each generated trade statistics
we evaluate irrationality indices $\omega_{G,b} =
\omega_G(\widetilde{TS}_b)$ and $\omega_{H,b} =
\omega_H(\widetilde{TS}_b)$. The power of GARP test is then
estimated as
\begin{gather*}
  \widehat{W}_G = \frac1B\sum_{b=1}^B\I(\omega_{G,b} > 1),
\end{gather*}
where $\I(\omega_{G,b} > 1)$ is equal to 1 if $\omega_{G,b} > 1$ and
equal to 0 otherwise. The power of HARP test is estimated by
\begin{gather*}
  \widehat{W}_H = \frac1B\sum_{b=1}^B\I(\omega_{H,b} > 1).
\end{gather*}

These estimates depend on the original trade statistics $\{(P^t,
X^t)\}_{t=1}^T$ as well as the distribution of $(\Tilde P^1, \ldots,
\Tilde P^T)$. In this paper we model $\Tilde P^t$ in the following
way. First, we fit auto-regression model $AR(r_i)$ for $z_i^t =
\log\left(\frac{P_i^t}{P_i^{t-1}}\right)$ for each good. We chose
$r_i$ from $\{0, 1, 2\}$ using Akaike information criterion (see,
for example, \citep{Hayashi2000}). We put the limit on possible
values of $r_i$ because the small number of time periods. The
estimated model for the price of good $i$ is
\begin{gather}
\begin{aligned}
  z_i^t & = \hat\beta_{i,0} + \sum_{\tau = 1}^{\hat
  r_i}\hat\beta_{i,\tau}z_i^{t-\tau} + \varepsilon_{i,t}, \quad t
  \in\{\hat r_i + 2,\ldots, T\}, \\
  z_i^t & = \log\left(\frac{P_i^t}{P_i^{t-1}}\right), \quad t \in
  \{1,\ldots,\hat r_i + 1\},
\end{aligned}\label{AR}
\end{gather}
where $\varepsilon_{i,t}\sim N(0, \hat\sigma_i^2)$ are independent
random variables. We put hats on the parameters of $AR$ model as
well as its order to emphasize that these are estimates. The
distributions of $\Tilde P^t$ are such that $\tilde z_i^t =
\log\left(\frac{\Tilde P_i^t}{\Tilde P_i^{t-1}}\right)$ follow
(\ref{AR}).

We took $B = 20000$ and estimated $\widehat{W}_G$ and
$\widehat{W}_H$ for the groups selected by Hungarian trade
experts\footnote{The groups in table \ref{HungGroups} represent the
classification of goods; trade experts considered not only them.}
which satisfy HARP. These groups are shown in table
\ref{ForecastTable} in the next section. The values of
$\widehat{W}_H$ lie in $[0.97540, 0.99985]$. The values of
$\widehat{W}_G$ lie in $[0, 0.53080]$. The histogram for
$\widehat{W}_G$ is shown in figure \ref{GARPPower}. We see that the
power of GARP test is very low compared to that of HARP test. The
high power of HARP test implies that the fact that a certain group
of goods satisfies HARP is a valuable finding, because it implies a
strong relation between consumption and price. The low power of GARP
test implies that the GARP often holds true even for group of goods
with no relation between price and consumption. That is why GARP
should not be used as a tool for studying the structure of consumer
demand.

\begin{figure}[!h]
  \begin{center}
    \includegraphics[width = \linewidth]{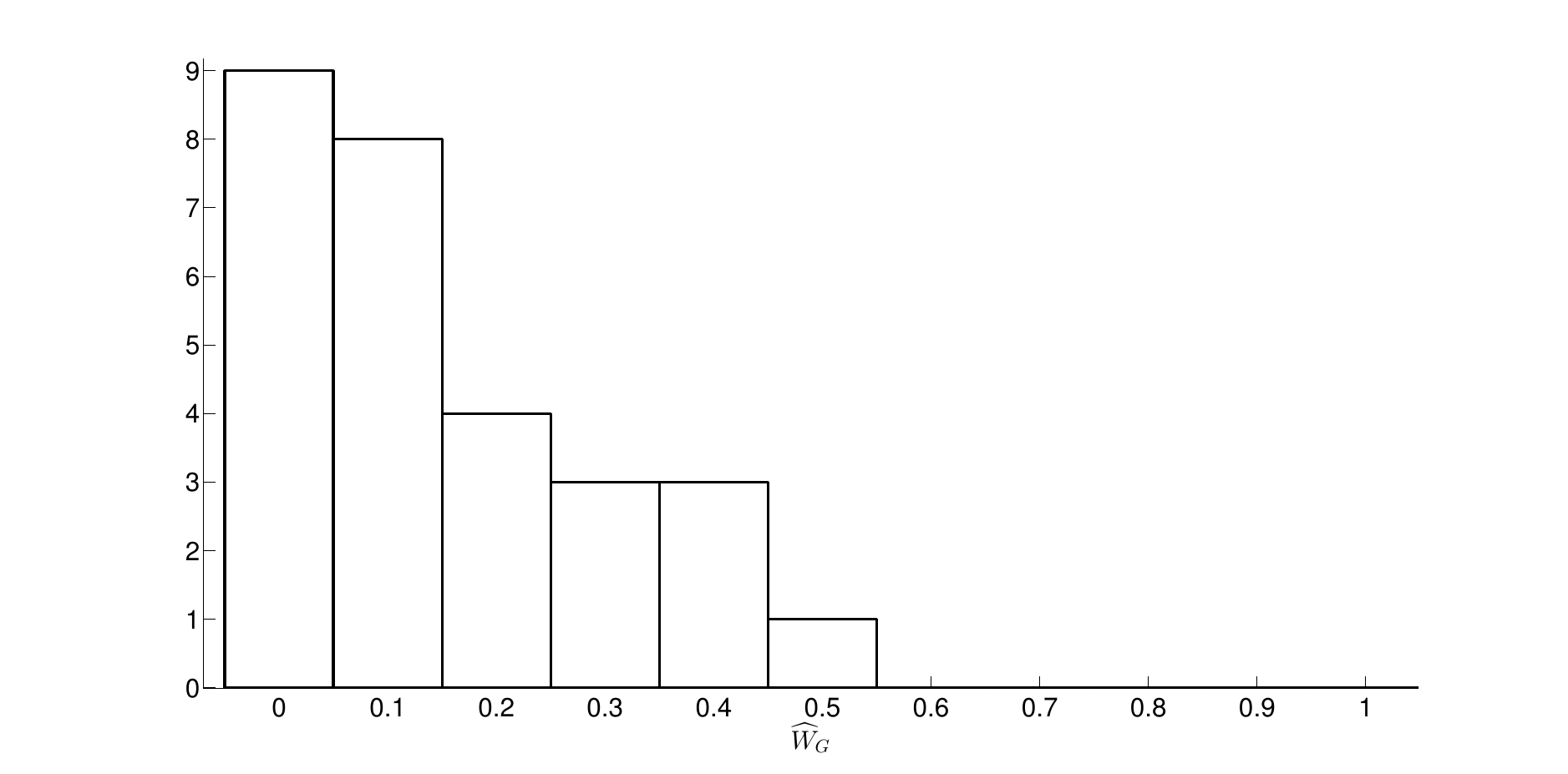}
    \caption{The histogram of $\widehat{W}_G$.}
    \label{GARPPower}
  \end{center}
\end{figure}

\section{Size of forecasting set}\label{SSize}

In this section we use the trade statistics of Hungary in order to
assess the size of the forecasting set for the last period price
given the actual value of last period consumption. The approach used
in this section is similar to that used in the previous section. The
groups are also the same.

Consider random trade statistics $\widetilde{TS} = \{(\Tilde P^t,
X^t)\}_{t=1}^T$, where $\Tilde P^t = P^t$ for $t
\in\{1,\ldots,T-1\}$. The last period price is a random vector with
uniformly distributed on $\{P\in\R^m_+ \mid \|P\| = 1\}$. The size
of forecasting set is measured in the following way. We generate $B$
realizations of $\widehat{TS}$: $\{\widetilde{TS}_b\}_{b=1}^B$. The
size of forecasting set based on GARP is measured by
\begin{gather*}
  \widehat{F}_G = \frac1B\sum_{b=1}^B\I(\widetilde{TS}_b \text{ satisfies
  GARP}),
\end{gather*}
where $\I(\widetilde{TS}_b \text{ satisfies GARP})$ is equal to 1 if
$\widetilde{TS}_b$ satisfies GARP and is equal to 0 otherwise. The
size of forecasting set based on HARP is measured similarly:
\begin{gather*}
  \widehat{F}_H = \frac1B\sum_{b=1}^B\I(\widetilde{TS}_b \text{ satisfies
  HARP}).
\end{gather*}

These measures are consistent estimated for the probability of
having a random price vector $P$ in the forecasting set based on
GARP and HARP respectively. The values of any of these measures
close to 1 imply that the corresponding forecasting set is close to
the set of all possible price vectors.

The estimated values of $\widehat{F}_G$ and $\widehat{F}_H$ for $B =
100000$ are shown in table \ref{ForecastTable}. We see that for the
majority of the considered groups the value of $\widehat{F}_G$ is
equal to one. This shows that the forecasting set based on GARP
often gives almost no information about future prices. The situation
is different for the forecasting set based on HARP. We see that this
forecasting set indeed provide some information on future prices as
it puts certain restrictions on them.

\begin{table}[!h]
  \begin{center}
    \caption{The sizes of forecasting sets for selected groups from Hungarian trade statistics.}
    \label{ForecastTable}
    \begin{tabular}{L{7cm}cc}
      \toprule
      Group description & $\widehat{F}_H$ & $\widehat{F}_G$ \\
      \midrule
      Meat $+$ meat products & 0.29239 & 0.76349 \\
      Same $+$ fish products & 0.2763 & 0.71422 \\
      Foodstuff & 0.27244 & 1 \\
      Same $+$ sweets $+$ fruits and vegetables $+$ beverages $+$
      tobacco & 0.15899 & 1 \\
      (1-67) + clothes + housing service + energy in home  +
      furnishings and home equipment & 0.1535 & 1 \\
      Same + health and hygiene & 0.15713 & 1 \\
      Same + transport and information & 0.15538 & 1 \\
      Same + education, culture, sport and recreation
      & 0.15969 & 1 \\
      All goods & 0.16237 & 1 \\
      \bottomrule
    \end{tabular}
  \end{center}
\end{table}

\section{Conclusion}\label{SConcl}

In this paper we have demonstrated an approach to construction of
economic indices based of homothetic axiom of revealed preference.
We have provided examples of computing these indices and building
the tree of economic indices which serves as a visualization of the
structure of consumer demand. We demonstrated that testing data for
consistency with this axiom is a better approach to studying the
structure of consumer demand that based on generalized axiom of
revealed preference. We provided two arguments for this.

The first argument is that there is high probability of having a
random group satisfying GARP. This means that the analysis based on
GARP reveals that almost any group is full in terms of substitution
and complementarity. Such outcome is useless and contradicts
intuition. Most probably by choosing absolutely random group we
should reveal that it does not represent an independent market
segment. The analysis based on HARP allows one to choose full groups
of not very great number of goods more accurately which makes it
proper tool for studying the market segmentation.

The second argument is that the power of GARP test is very low
compared to that of HARP test. We compared the power of HARP and
GARP test on the groups which satisfy HARP. The fact that these
groups satisfy HARP provides additional evidence in favor of their
completeness. However, the power of GARP test on these groups is
very low. This means that GARP test does not catch strong relation
between prices and consumption of the goods in the groups despite we
have strong evidence in favor of existence of such relation.

We also have demonstrated that HARP test is suitable for forecasting
future prices. The resulting forecasting set is not a trivial one.
This is not the case for GARP test. The forecasting set based on
GARP is close to a trivial one. Moreover, in many cases it coincides
with it.

\section*{Appendix 1}

\textbf{Proof of corollary \ref{HARPInfCol1}.} The relations
(\ref{HARPIndIneq1}) and (\ref{HARPIndIneq2}) imply that for all $X
\in \inter\R^m_+$ and $P \in \R^m_+$ such that $Q(P) > 0$
\begin{gather*}
  \frac1{Q(P(X))}\dpr{P(X)}{X} \leqslant \frac1{Q(P)}\dpr{P}{X}.
\end{gather*}
The division by $Q(P(X))$ is correct because if $X \in \inter\R_+^m$
then $F(X) > 0$n and $\dpr{P(X)}{X} > 0$. Therefore,
(\ref{HARPIndIneq2}) implies $Q(P(X)) > 0$. Substituting $P$ with
$P(Y)$ leads to
\begin{gather*}
  \frac1{Q(P(X))}\dpr{P(X)}{X} \leqslant \frac1{Q(P(Y))}\dpr{P(Y)}{X}
  \quad \forall X, Y \in \inter\R_+^m.
\end{gather*}
Therefore, $\lambda(X) = \frac1{Q(P(X))}$ satisfies
(\ref{HARPSysInf}). Since $Q(P(X))$ is continuous so is
$\lambda(X)$.

\begin{flushright}
  \textit{Q.E.D.}
\end{flushright}

\textbf{Proof of corollary \ref{HARPInfCol2}.} Define $F(X)$ and
$Q(P)$ by
\begin{gather*}
  F(X) = \inf\left\{\left.\lambda(Y)\dpr{P(Y)}{X} \right| Y \in
  \R_+^m\right\},
  Q(P) = \inf\left\{\left.\frac{\dpr{P}{Y}}{F(Y)}\right|Y\geqslant0, F(Y) >
  0\right\},
\end{gather*}
where $\lambda(Y)$ satisfies (\ref{HARPSysInf}). Then the functions
$F(X)$ and $Q(P)$ satisfy (\ref{HARPIndIneq1}). Since $\lambda(Y)$
satisfies (\ref{HARPSysInf}),
\begin{gather*}
  F(X) = \lambda(X)\dpr{P(X)}{X},
\end{gather*}
and
\begin{multline*}
  Q(P(X)) = \inf\left\{\left.\frac{\dpr{P(X)}{Y}}{\lambda(Y)\dpr{P(Y)}{Y}}\right|Y\geqslant0, \lambda(Y)\dpr{P(Y)}{Y} >
  0\right\} \\
  = \frac1{\lambda(X)} \inf\left\{\left.\frac{\lambda(X)\dpr{P(X)}{Y}}{\lambda(Y)\dpr{P(Y)}{Y}}\right|Y\geqslant0, \lambda(Y)\dpr{P(Y)}{Y} >
  0\right\} = \frac1{\lambda(X)}.
\end{multline*}

\begin{flushright}
  \textit{Q.E.D.}
\end{flushright}

\textbf{Proof of corollary \ref{VSCoroll}.} We show that functions
$F(X)$ and $Q(P)$ defined as
\begin{gather*}
  F(X) = F_0(F_1(X_1),\ldots,F_K(X_K), X_{K+1}), \\
  Q(P) = Q_0(Q_1(P_1),\ldots,Q_K(P_K), P_{K+1})
\end{gather*}
are from $\Phi_H$ and satisfy (\ref{HARPIndIneq1}) and
(\ref{HARPIndIneq2}). The fact that these functions are from
$\Phi_H$ is implied from the definition of functions $Q_k$ and $F_k$
($k=1,\ldots,K+1$). The inequality (\ref{HARPIndIneq1}) follows from
(\ref{VSCoroll1}) and (\ref{VSCoroll3}):
\begin{multline*}
  Q(P)F(X) = Q_0(Q_1(P_1),\ldots,Q_K(P_K),P_{K+1})F_0(F_1(X_1),\ldots,F_K(X_K),
  X_{K+1}) \\
  \leqslant \sum_{k=1}^KQ_k(P_k)F_k(X_k) + \dpr{P_{K+1}}{X_{K+1}}
  \leqslant \dpr{P}{X}
\end{multline*}
The identity (\ref{HARPIndIneq2}) follows from (\ref{VSCoroll4}).
Theorem~\ref{HARPThInf} implies that the inverse demand functions
$P(X)$ are rationalizable in $\Phi_H$. The identities
(\ref{VSCoroll5}) and (\ref{VSCoroll6}) are satisfied by definition
of functions $F(X)$ and $Q(P)$.

\begin{flushright}
  \textit{Q.E.D.}
\end{flushright}

\textbf{Proof of theorem \ref{GGARPTh}}. 1) $\Rightarrow$ 2). This
part represent a modified version of the part of the proof of
theorem~\ref{AVGARPTh} given in \citep{Varian1982a}. Let $\max(I)$
be the index of maximal element with respect to the binary relation
$\mathit{R}^*(\omega)$ which is the transitive closure of binary
relation $\mathit{R}(\omega)$ defined on
$\{X^t\}_{t=1}^T\times\{X^t\}_{t=1}^T$ as
\begin{gather*}
  X^t\mathit{R}(\omega)X^{\tau} \Leftrightarrow \dpr{P^t}{X^t}
  \geqslant \omega\dpr{P^t}{X^{\tau}}.
\end{gather*}
In other words,
\begin{gather*}
  \forall t\in I \quad X^t\mathit{R}^*(\omega)X^{\max(I)} \Rightarrow
  X^{\max(I)}\mathit{R}^*(\omega)X^t.
\end{gather*}

Consider the following algorithm.

Output: A set of numbers $U^t$, $\lambda^t > 0$ ($t =
\overline{1,T}$).
\begin{itemize}
  \item[1)] $I = \{1,\ldots,T\}$, $B=\varnothing$.
  \item[2)] Let $m = \max(I)$.
  \item[3)] Set $E = \{t \in I \mid X^t\RP X^m\}$. If $B = \varnothing$,
  set $U^m = \lambda^m = 1$ and go to 6). Otherwise go to 4).
  \item[4)] Set $U^m = \min\limits_{t\in E}\min\limits_{\tau\in
  B}\min\left\{U^{\tau} + \lambda^{\tau}\left(\omega\dpr{P^{\tau}}{X^t} -
  \dpr{P^{\tau}}{X^{\tau}}\right), U^{\tau}\right\}$.
  \item[5)] Set $\lambda^m = \max\limits_{t\in E}\max\limits_{\tau\in
  B}\max\left\{\dfrac{U^{\tau} - U^m}{\omega\dpr{P^t}{X^{\tau}} -
  \dpr{P^t}{X^t}}, 1\right\}$.
  \item[6)] Set $U^t = U^m$, $\lambda^t = \lambda^m$ for all $t \in
  E$.
  \item[7)] Set $I = I \setminus E$, $B = B \cup E$. If $I =
  \varnothing$, stop. Otherwise, go to 2).
\end{itemize}

Let us prove that if trade statistics satisfies GARP($\omega$), then
this algorithm provides the solution to (\ref{GGARPSys}). The
algorithm is an iterative process. We show that after each iteration
of step 6) the constructed $U$'s and $\lambda$'s satisfy
corresponding inequalities of (\ref{GGARPSys}). Namely, we show that
\begin{gather*}
  (a) \quad U^t \leqslant U^{\tau} +
  \lambda^{\tau}\left(\omega\dpr{P^{\tau}}{X^t} -
  \dpr{P^{\tau}}{X^{\tau}}\right) \quad \forall \tau\in B, \, t\in E,
  \\
  (b) \quad U^{\tau} \leqslant U^t +
  \lambda^t\left(\omega\dpr{P^t}{X^{\tau}} - \dpr{P^t}{X^t}\right)
  \quad \forall \tau\in B, \, t\in E, \\
  (c) \quad U^t \leqslant U^{\tau} +
  \lambda^{\tau}\left(\omega\dpr{P^{\tau}}{X^t} -
  \dpr{P^{\tau}}{X^{\tau}}\right) \quad \forall t,\tau \in E,\, t\ne\tau.
\end{gather*}

\textit{Proof of (a)}: By step 4) of the algorithm:
\begin{gather*}
  U^t = U^m \leqslant U^{\tau} +
  \lambda^{\tau}\left(\omega\dpr{P^{\tau}}{P^t} -
  \dpr{P^{\tau}}{X^{\tau}}\right) \quad \forall \tau\in B, \, t\in E.
\end{gather*}

\textit{Proof of (b)}: For this we need to use step 5). Notice that
$\omega\dpr{P^t}{X^{\tau}} > \dpr{P^t}{X^t}$ for all $\tau \in B$.
If not, then $X^t\RP X^{\tau}$ for some $\tau\in B$. But then $t$
would have been moved into $B$ before $\tau$ was and we have
contradiction with $t \in E$.

Hence, the division is well defined and
\begin{gather*}
  \lambda^t = \lambda^m \geqslant \frac{U^{\tau} -
  U^t}{\omega\dpr{P^t}{X^{\tau}} - \dpr{P^t}{X^t}} \quad \forall \tau\in B, \, t\in E.
\end{gather*}

Cross multiplying:
\begin{gather*}
  \lambda^t\left(\omega\dpr{P^t}{X^{\tau}} - \dpr{P^t}{X^t}\right)
  \geqslant U^{\tau} - U^t \quad \forall \tau\in B, \, t\in E,
\end{gather*}
which proves (b).

\textit{Proof of (c)}: Note that if $t, \tau \in E$, then
$\omega\dpr{P^{\tau}}{X^t} \geqslant \dpr{P^{\tau}}{X^{\tau}}$. If
not, then $\dpr{P^{\tau}}{X^{\tau}} > \omega\dpr{P^{\tau}}{X^t}$,
which contradicts GARP($\omega$). Indeed, $t \in E$ implies
$X^{t}\RP X^m$, $\tau \in E$ implies $X^{\tau} \RP X^m$ which (by
definition of $m$) implies $X^m \RP X^{\tau}$, so $X^{t} \RP
X^{\tau}$ and GARP($\omega)$ implies $\dpr{P^{\tau}}{X^{\tau}}
\leqslant \omega\dpr{P^t}{X^{\tau}}$. Now for all $t,\tau \in E$:
\begin{gather*}
  U^t = U^{\tau}, \quad \lambda^{\tau} = \lambda^m > 0,
\end{gather*}
so
\begin{gather*}
  U^t \leqslant U^{\tau} +
  \lambda^{\tau}\left(\omega\dpr{P^{\tau}}{X^t} -
  \dpr{P^{\tau}}{X^{\tau}}\right).
\end{gather*}

2) $\Rightarrow$ 1). Let $\{U^t, \lambda^t\}_{t=1}^T$ with
$\lambda^t > 0, \forall t \in \{1,\ldots, T\}$ satisfy
(\ref{GGARPSys}), and $\dpr{P^t}{X^t} \geqslant \omega
\dpr{P^t}{X^{t_1}}$, $\dpr{P^{t_1}}{X^{t_1}} \geqslant \omega
\dpr{P^{t_1}}{X^{t_2}}$, $\ldots$, $\dpr{P^{t_k}}{X^{t_k}} \geqslant
\omega \dpr{P^{t_k}}{X^s}$. Taking into account that $\lambda^t > 0$
$\forall t \in \{1,\ldots,T\}$, these inequalities imply that $U^t
\leqslant U^s$. This implies $\dpr{P^s}{X^s} \leqslant \omega
\dpr{P^s}{X^t}$ which completes the proof of this part.

\begin{flushright}
  \textit{Q.E.D.}
\end{flushright}

\textbf{Proof of theorem \ref{GHARPandGGARP}.} 1) $\Rightarrow$ 2).
HARP($\omega$) is invariant with respect to change of scales.
Therefore if the trade statistics $\{(P^t,X^t)\}_{t=1}^T$ satisfies
HARP($\omega$) then the trade statistics
$\{(P^t,\mu^tX^t)\}_{t=1}^T$ also satisfies HARP($\omega$) for any
$\mu^t$ $(t\in\{1,\ldots,T\})$. This implies that
$\{(P^t,\mu^tX^t)\}_{t=1}^T$ satisfies GARP($\omega$).

2) $\Rightarrow$ 1). Fix some subset of indices $(t_1,\ldots,t_k)
\subset \{1,\ldots,T\}$ such that there are no two identical
indices. Select $\mu^{t_1}$, $\ldots$, $\mu^{t_k}$ so that
\begin{gather*}
  \dpr{P^{t_1}}{\mu^{t_1}X^{t_1}} =
  \omega\dpr{P^{t_1}}{\mu^{t_2}X^{t_2}}, \\
  \dpr{P^{t_2}}{\mu^{t_2}X^{t_2}} =
  \omega\dpr{P^{t_2}}{\mu^{t_3}X^{t_3}}, \\
  \ldots \\
  \dpr{P^{t_{k-1}}}{\mu^{t_{k-1}}X^{t_{k-1}}} =
  \omega\dpr{P^{t_{k-1}}}{\mu^{t_k}X^{t_k}}.
\end{gather*}
It is possible to do that because there are $k - 1$ equalities to
satisfy by choosing $k$ numbers. Then GARP($\omega$) implies
\begin{gather*}
  \dpr{P^{t_k}}{\mu^{t_k}X^{t_k}} \leqslant
  \omega\dpr{P^{t_k}}{\mu^{t_1}X^{t_1}}.
\end{gather*}
Therefore,
\begin{gather*}
  \frac{\dpr{P^{t_1}}{X^{t_2}}}{\dpr{P^{t_1}}{X^{t_1}}}\frac{\dpr{P^{t_2}}{X^{t_3}}}{\dpr{P^{t_2}}{X^{t_2}}}
  \ldots \frac{\dpr{P^{t_k}}{X^{t_1}}}{\dpr{P^{t_k}}{X^{t_k}}} \geqslant \frac1{\omega^k}.
\end{gather*}

\begin{flushright}
  \textit{Q.E.D.}
\end{flushright}

\section*{Appendix 2}

Consider a group of three goods the demand on which is given by the
three Engel curves
\begin{gather*}
  q_1(x) = (1,\,\ve,\,\ve)'x, \\
  q_2(x) = (\ve,\,1,\,\ve)'x, \\
  q_3(x) = (\ve,\,\ve,\,1)'x,
\end{gather*}
which correspond to the following price vectors:
\begin{gather*}
  P^1 = \left(2\,1\,4\right)', \\
  P^2 = \left(2\,1\,2\right)', \\
  P^3 = \left(2\,2\,1\right)'.
\end{gather*}

Observed demands are given by $X^t(\ve) = q_t(1)$. The matrix
$PX(\ve)$ with elements $px_{\tau t}(\ve) = \dpr{P^\tau}{X^t(\ve)}$
is equal to
\begin{gather*}
  PX(\ve) = \left(
  \begin{array}{ccc}
    2+5\ve & 1+6\ve & 4+3\ve \\
    2+3\ve & 1+4\ve & 2+3\ve \\
    2+3\ve & 2+3\ve & 1+4\ve
  \end{array}
  \right).
\end{gather*}

If $\ve = 0$ then the direct revealed preference relation is given
by
\begin{gather*}
  R(0) = \left(
  \begin{array}{ccc}
    1 & 1 & 0 \\
    0 & 1 & 0 \\
    0 & 0 & 1
  \end{array}
  \right),
\end{gather*}
with transitive closure $R^*(0) = R(0)$. The trade statistics
$\{(P^t, X^t(0))\}_{t=1}^T$ satisfies GARP and HARP. In order to
simplify the calculations we set $\ve = 0$. Having zero demand may
seem unnatural, however the forecasting sets for small enough
positive $\ve$ should differ only a little unless setting $\ve>0$
leads to failure of GARP\footnote{This is not a problem for HARP
since this axiom is a fixed set of non-strict inequalities.}. It may
be shown that $R(\ve) = R(0)$ for $\ve < 1$. Therefore, choosing
$\ve$ from $(0,1)$ does not lead to failure of GARP.

Let $P^4 = (1\,1\,1)'$ and $x_4 = 2$. Since Engel curves are rays,
$K_H^1(P^{4}, x_{4}) = \widetilde{K}_G^1(P^4, x_4)$. Therefore, if
$\widetilde{K}_G^1(P^4, x_4) = G(P^{4}, x_{4}, \{\tilde
x_t\}_{t=1}^3)$, then $G(P^{4}, x_{4}, \{\tilde x_t\}_{t=1}^3) =
K_H^1(P^{4}, x_{4})$. We show that this is not true.

We start with $G(P^{4}, x_{4}, \{\tilde x_t\}_{t=1}^3)$. The
intersection demands are given by
\begin{gather*}
  \dpr{P^4}{q_t(\tilde x_t)} = \tilde x_t = 2.
\end{gather*}
and the set $\{X\in\R_+^m \mid X \in G(P^{4}, x_{4}, \{\tilde
x_t\}_{t=1}^3)$ is given by
\begin{align*}
  & X_1 \leqslant 2 - \frac32 X_2 \\
  & X_3 = 2 - X_1 - X_2 \\
  & X_1 \geqslant 0, X_2 \geqslant 0.
\end{align*}

The set $\{X\in\R_+^m \mid X \in K_H^1(P^{4}, x_{4})$ is given by
\begin{align*}
  & X_3 = 2 - X_1, \\
  & X_2 = 0, \\
  & X_1 \geqslant 0,\, X_1 \leqslant 2.
\end{align*}

We see that the set $K_H^1(P^4, x_4)$ is a proper subset of
$G(P^{4}, x_{4}, \{\tilde x_t\}_{t=1}^3)$. This implies that
(\ref{BBCTheorem}) is not true. If we set $\ve = 0.5$ then $G(P^{4},
x_{4}, \{\tilde x_t\}_{t=1}^3)$ is been given by
\begin{align*}
  & X_1 \leqslant 1.75 - 1.5X_2 \\
  & X_2 \leqslant 1 \\
  & X_1 + X_2 \geqslant 1 \\
  & X_3 = 2 - X_1 - X_2 \\
  & X_1 \geqslant 0,\, X_2 \geqslant 0,
\end{align*}
while the set $K_H^1(P^4, x_4)$ is given by
\begin{align*}
  & X_1 \leqslant 1.75 - 1.5X_2 \\
  & X_2 \leqslant 0.625 \\
  & X_1 + X_2 \geqslant 1 \\
  & X_3 = 2 - X_1 - X_2 \\
  & X_1 \geqslant 0,\, X_2 \geqslant 0.
\end{align*}
The set $K_H^1(P^4, x_4)$ is still a proper subset of $G(P^{4},
x_{4}, \{\tilde x_t\}_{t=1}^3)$.

\bibliographystyle{apalike}
\bibliography{MyBiblio}

\end{document}